\RequirePackage{etoolbox}
\csdef{input@path}{{style/}{graphics/}}
\documentclass[MSNbibl,seceqn,number,citesort,dvips]{arxstspdf}
\usepackage{multirow}
\usepackage{graphicx}

\usepackage{url,breakurl}
\usepackage{flushend}
\usepackage{stfloats}

\volume{30}
\issue{1}
\pubyear{2015}
\firstpage{1}
\lastpage{25}
\doi{10.1214/14-STS506}
\docsubty{FLA}

\makeatletter
\newproclaim{remark}{Remark}

\newtheorem{lemma}{Lemma}[section] 

\newproclaim{definition}{Definition}[section]

\newcommand{\eps}{\varepsilon}
\newcommand{\bitem}{\begin{itemize}}
\newcommand{\eitem}{\end{itemize}}
\newcommand{\goto}{\rightarrow}
\newcommand{\margmax}{\mathop{\operatorname{argmax}}}
\newcommand{\that}{\hat{t}}

\makeatother

\begin{document}
\begin{frontmatter}

\dochead{Special Invited Paper}
\title{Higher Criticism for Large-Scale Inference, Especially for Rare
and Weak Effects}
\runtitle{Higher Criticism for Rare and Weak Effects}

\begin{aug}
\author[A]{\fnms{David} \snm{Donoho}\ead[label=e1]{donoho@stanford.edu}}
\and
\author[B]{\fnms{Jiashun} \snm{Jin}\corref{}\ead[label=e2]{jiahsun@stat.cmu.edu}}

\runauthor{D. Donoho and J. Jin}

\address[A]{David Donoho is
Professor of Statistics,
Statistics Department,
Stanford University,
Stanford, California 94305,
USA \printead{e1}.}
\address[B]{Jiashun Jin is Professor of Statistics,
Statistics Department,
Carnegie Mellon University,
5000 Forbes Avenue,
Pittsburgh, Pennsylvania 15213,
USA \printead{e2}.}
%
%
\affiliation{Stanford University and Carnegie Mellon University}
\end{aug}

%
\begin{abstract}
In modern high-throughput data analysis, researchers perform a large
number of
statistical tests, expecting to find perhaps a small fraction of
significant effects
against a predominantly null background. Higher Criticism (HC) was
introduced to
determine whether there are \textit{any} nonzero effects; more recently,
it was applied to
feature selection, where it provides a method for selecting useful
predictive features from
a large body of potentially useful features, among which only a rare
few will prove
truly useful.

In this article, we review the basics of HC in both the testing and
feature selection settings.
HC is a flexible idea, which adapts easily to new situations; we point
out simple adaptions to clique detection
and bivariate outlier detection. HC, although still early in its development,
is seeing increasing interest from practitioners; we illustrate this
with worked examples.
HC is computationally effective, which gives it a nice leverage in the
increasingly more relevant ``Big Data'' settings we see today.

We also review the underlying theoretical ``ideology'' behind HC.
The \textit{Rare/Weak} (RW) model is a theoretical framework
simultaneously controlling the size and prevalence of
useful/significant items among the useless/null bulk.
The RW model shows that
HC has important advantages over better known procedures
such as False Discovery Rate (FDR) control and Family-wise Error
control (FwER),
in particular, certain optimality properties. We discuss the rare/weak
\textit{phase diagram},
a way to visualize clearly the class of RW settings where the true
signals are so
rare or so weak that detection and feature selection are simply impossible,
and a way to understand the known optimality properties of HC.
\end{abstract}
%
%
\end{frontmatter}

\section{Introduction} \label{sec:Intro}

A data deluge is now flooding scientific and technical
work \cite{datadeluge}.
In field after field, high-throughput devices
gather many measurements per individual;
depending on the field, these could be gene expression levels,
or spectrum levels, or peak detectors or wavelet transform
coefficients; there could be
thousands or even millions of different feature measurements per single subject.

High-throughput measurement technology automatically measures
systematically generated features and contrasts; these features are not
custom-designed for any one project.
Only a small proportion of the
measured features are expected to be relevant for the research in question,
but researchers do not know in advance which those will be;
they instead measure every contrast fitting within their systematic scheme,
intending later to identify a small fraction of relevant ones post-facto.

This flood of high-throughput measurements is
driving a new branch of statistical practice:
what Efron \cite{EfronLSI} calls \textit{Large-Scale Inference} (LSI).
For this paper, two specific LSI problems are of interest:
\bitem
\item\textit{Massive multiple testing for sparse intergroup differences.}
Here we have two groups, a treatment and a control,
and for each measured variable we test whether the two groups are
different on that measurement, obtaining, say, a $P$-value per feature.
Of course, many individual features are unrelated to the specific
intervention being studied,
and those would be expected to show no significant differences---but
we do not know which
these are.
We expect that even when there are true inter-group differences,
only a small fraction of measured features will be affected---but,
again, we do not know which features they are.
We therefore use the whole collection of $P$-values to correctly decide
if there is any difference between the two groups.
\item\textit{Sparse feature selection.} A large number of features are
available for
training a linear classifier, but we expect that most of those
features are in fact useless for
separating the underlying classes. We must decide which features
to use in designing a class prediction rule.
\eitem
Higher Criticism (HC) and its elaborations can be useful in both of
these LSI settings;
under a particular asymptotic model discussed below---the \textit{Asymptotic Rare/Weak} (ARW)
\textit{model}---HC
offers theoretical optimality in selecting features.
In this paper we will review the basic notions of
HC, some variations and settings where it applies.
HC is a flexible idea and can be adapted to a range of
new problem areas; we briefly discuss three simple examples.

\subsection{HC Basics}
\label{subsec:HCbasic}
John Tukey \cite{Tukey76,Tukey89,Tukey94} coined the term ``Higher
Criticism''\footnote{In mid-twentieth century humanities
studies, the term Higher Criticism became popular to label a certain
school of Biblical scholarship.}
and motivated it by the following story. A~young scientist administers
a total of $250$ independent tests,
out of which $11$ are significant at the level of $5\%$. The youngster
is excited about the findings and plans to trumpet them
until a senior researcher tells him that, even
in the purely null case, one would expect to have $12.5$
significances. In that sense, finding only $11$ significances is
actually disappointing.
Tukey proposes a kind of
\textit{second-level significance testing}, based on the statistic
\begin{eqnarray*}
\mathrm{HC}_{N, 0.05} &=& \sqrt{N} (\mbox{Fraction significant at $5\%$} -
0.05)\\
&&{}/ \sqrt{0.05 \times0.95},
\end{eqnarray*}
where $N = 250$ is the total number of tests. Obviously this score has
an interpretation similar to
$Z$- and $t$- statistics,
so Tukey suggests that a value of $2$ or larger indicates \textit
{significance of the overall body of tests}.
In Tukey's example,
\[
\mathrm{HC}_{N, 0.05} = - 0.43.
\]
If the young researcher really ``had something'' this score should be
strongly positive,
for example, $2$ or more, but here the score is negative, implying that
the overall body of the evidence is consistent with the null
hypothesis of no difference.
Donoho and Jin \cite{DJ04} saw that in the modern context of
large $N$ and rare/weak signals, it was advantageous to generalize
beyond the single significance level $\alpha=0.05$.
They maximized over \textit{all} levels $\alpha$ between $0$ and some
preselected upper bound $\alpha_0 \in(0, 1)$.
So generalize Tukey's proposal and set
\begin{eqnarray*}
\mathrm{HC}_{N, \alpha} &= &\sqrt{N} (\mbox{Fraction significant at $\alpha$}
- \alpha)\\
&&{}/
\sqrt{\alpha\times(1 - \alpha)}.
\end{eqnarray*}
If the overall body of tests is significant, then we expect $\mathrm{HC}_{N,
\alpha}$ to be large for \textit{some} $\alpha$. Otherwise,
we expect $\mathrm{HC}_{N, \alpha}$ to be small over all $\alpha$ in a wide
range. In other words, the significance of
the overall body of test is captured in the following Higher Criticism
statistic:
%
\begin{equation}
\label{DefineHC1} \mathrm{HC}_N^* = \max_{\{ 0 \leq\alpha\leq\alpha_0 \}}
\mathrm{HC}_{N, \alpha},
\end{equation}
where $\alpha_0 \in(0,1)$ is a tuning parameter we often set at
$\alpha
_0 = 1/2$.

%
\begin{figure}[b]

\includegraphics{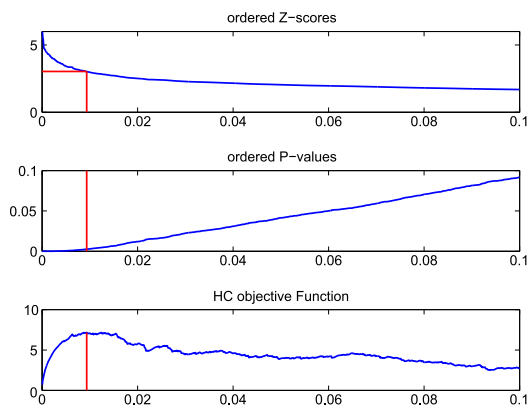}

\caption{Illustration of HC. The component score maximizing the HC
objective is located
at the red line. Bottom panel: the HC objective function $\mathrm{HC}_{N,i}$
versus $i/N$. Middle panel:
the underlying $P$-vales $\pi_{(i)}$ versus $i/N$. Top panel: the corresponding
ordered $Z$-scores $z_{(i)}$ versus $i/N$.}
\label{fig:IllustHC}
\end{figure}

Higher Criticism (HC) can be computed efficiently as follows.
Consider a total of $N$ uncorrelated tests:
\begin{itemize}
\item For $i = 1, 2, \ldots, N$, get the corresponding individual
$P$-values $\pi_i$,
producing in all a body of $P$-values $\pi_1, \pi_2, \ldots, \pi_N$.
\item Sort the $P$-values in the ascending order:
\[
\pi_{(1)} < \pi_{(2)} < \cdots< \pi_{(N)}.
\]
\item The Higher Criticism statistic in (\ref{DefineHC1}) can be
equivalently written as follows:
%
\begin{eqnarray}
\label{altHC} \mathrm{HC}_N^* &= &\max_{\{ 1 \leq i \leq\alpha_0 N\}}
\mathrm{HC}_{N, i},
\nonumber
\\[-8pt]
\\[-8pt]
\nonumber
 \mathrm{HC}_{N,
i} &\equiv&\sqrt{N} \frac{(i/N) - \pi_{(i)}}{\sqrt{\pi_{(i)} (1 -
\pi_{(i)})}}.
\end{eqnarray}
\end{itemize}
In words, we are looking at a test for equality of a binomial proportion
$\pi_{(i)}$ to an expected value $i/N$, maximizing this statistic across
a range of $i$. We think that the evidence against being purely null is
located somewhere in this range, but we cannot say in advance where that
might be.
The computational cost of $\mathrm{HC}$ is $O(N\log(N))$ and is moderate.

Figure~\ref{fig:IllustHC} illustrates the definition of Higher Criticism.
Consider an example where the (one-sided) $P$-values $\pi_i$ are
produced by $Z$-values $z_i$ through $\pi_i =1 - \Phi(z_i)$,
$1 \leq i \leq N$, where $\Phi$ denotes the CDF of $N(0,1)$. The first
panel shows the sorted $Z$-values in the descending order, the second
panel shows the sorted $P$-values, and the last panel shows $\mathrm{HC}_{N,
i}$. In this example, $\mathrm{HC}_N^* = 7.1$, reached by $\mathrm{HC}_{N, i}$ at $i =
0.0085 \times N$.

\begin{remark*} Asymptotic theory shows that the component scores
$\mathrm{HC}_{N,i}$ can be poorly behaved for $i$ very small (e.g., $1$ or $2$).
We often recommend the following modified version:
\[
\mathrm{HC}_N^+ = \max_{\{ 1 \leq i \leq\alpha_0 N\dvtx \pi_{(i)} > 1/N \}} \mathrm{HC}_{N, i}.
\]
\end{remark*}

\begin{remark*}
As the last remark illustrates,
small variations on the above prescription will
sometimes be useful, for example,
the modification of the underlying $Z$-like scores,
in (\ref{HCGF1})--(\ref{HCGF2}) below. Moreover, several other
statistics such as Berk--Jones and Average Likelihood Ratio offer
cognates or substitutes; see Section~\ref{sec:HCLike} below.
The real point of HC is less the specific definition (\ref{altHC})
and more a viewpoint
about the nature of evidence against the null hypothesis,
namely, that although the evidence may be cumulatively
substantial, it is diffuse, individually very
weak and affecting a relatively small fraction of the
individual $P$-values or $Z$-scores in our study.

So HC can be viewed as a family of methods for which the above definitions
give a convenient entry point. To make utterly clear, when needed
we label definition (\ref{altHC}) the Orthodox Higher Criticism (OHC).
\end{remark*}

\subsection{The Rare/Weak Effects Viewpoint and Phase Diagram}

Effect \textit{sparsity} was proposed as a useful hypothesis already in
the 1980s by
Box and Meyer \cite{Box}; it proposes that relatively few of the observational
units or factorial levels can be expected to show any difference from a global
null hypothesis of no effect, and that {a priori} we have no
opinion about which
units or levels those might be.

The Effect \textit{weakness} hypothesis assumes that
individual effects are not individually strong enough to be
detectable, once traditional
multiple comparisons ideas are taken into account.\footnote{For
example, Bonferroni-based
family-wise error rate control.}

The Rare/Weak viewpoint combines \textit{both} hypotheses in
analysis of large-scale experiments; it is intended to be a
flexible concept and to vary from one setting to another.

The next section operationalizes these ideas in a specific model,
where $N$ independent $Z$-scores follow a mixture with
a fraction $(1 - \eps)$ which are truly null effects and so distributed
$N(0,1)$,
while the remaining $\eps$ fraction
have a common effect size $\tau$ and are distributed $N(\tau,1)$. In
this situation,
the Rare/Weak viewpoint studies the regime where $\eps$ is small, the
locations of the
nonzero effects are scattered irregularly through the scores and the
effect size
$\tau$ is, at moderate $N$, only $2$ or $3$ standard deviations.

For large $N$ one can develop a precise theory;
see Section~\ref{sec:phase} below.
There we develop the \textit{Asymptotic Rare/Weak} (ARW) model,
a framework that assigns parameters to the \textit{rare} and \textit{weak}
attributes of a nonnull situation; a key phenomenon is that in the
two-dimensional parameter space there are three separate regions (or
\textit{phases})
where an inference goal is \textit{relatively easy}, \textit
{nontrivial but
possible}, and \textit{impossible} to achieve, correspondingly.
The ARW phase diagram offers revealing comparisons between
HC and other seemingly similar methods,
such as FDR control.

\section{HC for Detecting Sparse and Weak Effects} \label{sec:AoS}
In \cite{DJ04} Higher Criticism was originally proposed for detecting
sparse Gaussian mixtures.
Suppose we have $N$ test statistics $X_i$, $1 \leq i \leq N$
(reflecting many individual genes, voxels, etc.).
Suppose that these tests are standardized so that each individual test,
under its corresponding null
hypothesis, would have mean 0 and standard deviation 1.
We are interested in testing whether \textit{all} test statistics are
distributed
$N(0,1)$ versus the alternative that a small fraction is distributed as
normal with an elevated mean
$\tau$. In effect, we want an \textit{overall} test of a complete
null hypothesis:
%
\begin{equation}
\label{mixture1} H_0^{(N)}\dvtx\quad X_i \stackrel{
\mathrm{i.i.d.}} {\sim} N(0,1),\quad 1 \leq i \leq N,
\end{equation}
against an alternative in its complement,
%
\begin{eqnarray}
\label{mixture2} \quad H_1^{(N)}\dvtx \quad X_i \stackrel{
\mathrm{i.i.d.}} {\sim} (1 - \eps) N(0,1) + \eps N(\tau, 1),
\nonumber
\\[-8pt]
\\[-8pt]
 \eqntext{1 \leq i \leq N.}
\end{eqnarray}

%
\begin{figure}

\includegraphics{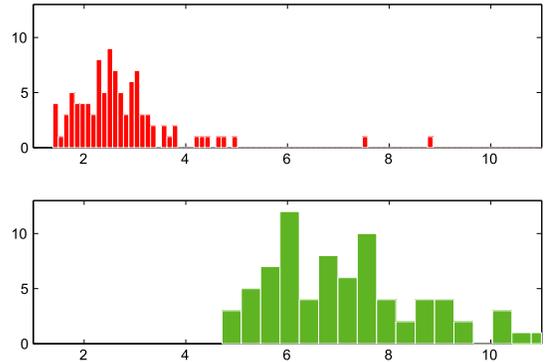}

\caption{Simulated Higher Criticism values. Top panel: simulation under
null hypothesis $H_0^{(N)}$.
Bottom panel: simulation under alternative hypothesis $H_1^{(N)}$.}
\label{fig:OAoHCHist}
\end{figure}

To use HC for such a case, we calculate the (one-sided) $P$-values by
\[
\pi_i = 1 - \Phi(X_i),\quad 1 \leq i \leq N.
\]
We then apply the basic definition of HC to the collection of
$P$-values.\footnote{If we
thought that under the alternative the mean might be either positive or
negative, we would of course use two-sided $P$-values.}

In Figure~\ref{fig:OAoHCHist}, we show the simulated HC values of
$H_0^{(N)}$ and $H_1^{(N)}$ based on $100$ independent
repetitions, where the parameters are set as $(N, \eps, \tau) = (10^6,
10^{-3}, 2)$. It is seen that the simulated HC values under $H_1^{(N)}$
are well separated
from those under $H_0^{(N)}$.

\subsection{Critical Value for Using HC as a Level-\texorpdfstring{$\alpha$}{alpha} Test}

Fix $0 < \alpha< 1$. To use HC as a level-$\alpha$ test, we must find
a critical value $h(N, \alpha)$ so that
\[
P_{H_0^{(N)}} \bigl\{ \mathrm{HC}_N^* > h(N, \alpha) \bigr\} \leq\alpha.
\]
$\mathrm{HC}^*_N$ can be connected with the maximum of a standardized
empirical process; see Donoho and Jin \cite{DJ04}. Using this
connection, it follows from \cite{Wellner86}, page 600, that as $N
\goto
\infty$, $b_N \mathrm{HC}_N^* - c_N$ and $b_N \mathrm{HC}_N^+ - c_N$ converge weakly to
the same limit---the standard Gumbel distribution, where $b_N = \sqrt{2
\log\log(N)}$ and $c_N = 2 \log\log(N) + (1/2) [\log\log\log(N)
- \log
(4\pi)]$.
As a result, for any fixed $\alpha\in(0,1)$ and $N \goto\infty$,
%
\begin{eqnarray}
\label{happrox} h(N, \alpha)& \approx& h_G(N,\alpha)
\nonumber
\\[-8pt]
\\[-8pt]
\nonumber
&=& \sqrt{2 \log
\log(N)} \bigl(1 + o(1)\bigr),
\end{eqnarray}
where $h_G(N,\alpha) = b_N^{-1}  [ c_N - \log\log(\frac{1}{1 -
\alpha}) ]$ (``G''  stands for Gumbel).
When $N$ is moderately large and $\alpha$ is moderately small, the
approximations may not be accurate enough,
and it is hard to derive an accurate closed-form approximation (even in
much simpler cases; see \cite{dasgupta} for nonasymptotic bound on
extreme values of normal samples). In such cases,
it is preferable to determine $h(N, \alpha)$ by simulations.
%
\begin{table*}
\tabcolsep=0pt
\caption{Simulated values $h(N, \alpha)$ based on $10^5$ repetitions.
Numbers in brackets are $h_G(N,\alpha)$}
\label{table:largen}
\begin{tabular*}{\textwidth}{@{\extracolsep{\fill}}lccccc@{}}
\hline
  & & \multicolumn{4}{c@{}} {$\bolds{N}$}
\\[-6pt]
& & \multicolumn{4}{c@{}} {\hrulefill}\\
{\textbf{Level}}& {\textbf{Statistic}}& $\bolds{10^3}$ &
 $\bolds{5 \times10^3}$ &$\bolds{2.5 \times10^4}$
& \multicolumn{1}{c@{}}{$\bolds{1.25 \times10^5}$}
\\
\hline
$\alpha= 0.05$ &$\mathrm{HC}_N^+$ & \phantom{0}3.17 (3.00) & \phantom{0}3.22 (3.08) & \phantom{0}3.26
(3.14) & \phantom{0}3.30 (3.19)
\\
&$\mathrm{HC}_N^*$ & \phantom{0}4.77 (3.00) & \phantom{0}4.73 (3.08) & \phantom{0}4.74 (3.14) &
\phantom{0}4.75 (3.19)
\\[3pt]
 $\alpha= 0.01$ &$\mathrm{HC}_N^+$ & \phantom{0}3.95 (3.83) & \phantom{0}3.97 (3.87) &
\phantom{0}3.96 (3.90) & \phantom{0}3.99 (3.93)
\\
&$\mathrm{HC}_N^*$ & 10.08 (3.83) & \phantom{0}9.88 (3.87) & 10.20 (3.90) &
\phantom{0}9.92 (3.93)
\\[3pt] $\alpha= 0.005$ &$\mathrm{HC}_N^+$ & \phantom{0}4.29 (4.18) & \phantom{0}4.28 (4.20) &
\phantom{0}4.26 (4.22) & \phantom{0}4.28 (4.24)
\\
 &$\mathrm{HC}_N^*$ & 13.78 (4.18) & 14.39 (4.20) &14.34 (4.22) &
13.95 (4.24)
\\[3pt]
$\alpha= 0.001$ &$\mathrm{HC}_N^+$ & \phantom{0}5.03 (5.00) & \phantom{0}5.02 (4.98)
&\phantom{0}4.98 (4.97) &\phantom{0} 4.98 (4.97)
\\
&$\mathrm{HC}_N^*$ & 30.27 (5.00) & 30.36 (5.02) &31.95 (4.97) &
31.49 (4.97)
\\
\hline
\end{tabular*}
\end{table*}

\begin{table*}[b]
\caption{MPM: malignant pleural mesothelioma. ADCA: adenocarcinoma.
ALL: acute lymphoblastic leukemia. AML: acute myelogenous leukemia}
\begin{tabular*}{\textwidth}{@{\extracolsep{\fill}}lccc@{}}
\hline
\textbf{Data name} & \textbf{\# training samples} &
\textbf{\# test samples} & \textbf{\# genes}
\\
\hline
Leukemia &27 (ALL), 11 (AML) & 20 (ALL), 14 (AML) & \phantom{0.}3571
\\
Lung cancer & 16 (MPM), 16 (ADCA) &15 (MPM), 134 (ADCA) & 12,533
\\
\hline
\end{tabular*}
\end{table*}

Table~\ref{table:largen} displays $h_G(N, \alpha)$ and $h(N,\alpha)$
[where $\alpha_0 = 1/2$ as in (\ref{DefineHC1})]
computed from $10^5$ independent simulations. One sees
that: (a) $h_G(N, \alpha)$ approximate the percentiles of $\mathrm{HC}_N^*$
poorly, but approximate those of $\mathrm{HC}_N^+$ reasonably well, especially
when $N$ get larger and $\alpha$ get smaller; (b) the tail of $\mathrm{HC}_N^*$
is fat but that of $\mathrm{HC}_N^+$ is relatively thin; (c) the percentiles of
$\mathrm{HC}_N^+$ and $\mathrm{HC}_N^*$ increase with $N$ only very slowly,
therefore, the values of $h(N,\alpha)$ for a few selected
$N$ represent those of a wide range of $N$. Very recently, Li and
Siegmund \cite{LiSiegmund} proposed a new
approximation to $h(N, \alpha)$ which is more accurate when $N$ is
moderately large.

\subsection{Two Gene Microarray Data Sets}
In Sections~\ref{subsec:realexample1} and \ref{sec:application}, we use
two standard
gene microarray data sets to help illustrate the use of HC: the lung
cancer data analyzed by Gordon {et al.} \cite{lungcancer},
and the leukemia data analyzed by Golub {et al.} \cite{leukemia}
(for the latter, we use the cleaned version published by Dettling \cite
{Dettling}, which contains measurements for $3571$ genes).
Both data sets are available at
\surl{http://www.stat.cmu.edu/\\\textasciitilde jiashun/Research/software/}. See Table~2,
where the partition of samples into the training set and the test set
is the same as in \cite{lungcancer} and~\cite{leukemia}, respectively.

\subsection{Detecting Rare and Weak Effects in Genomics and Genetics}
\label{subsec:genemicro}
When the genomics revolution began 10--15 years ago,
many scientists were hopeful that the common disease-common variant
hypothesis \cite{CDCVH} would apply. Under this hypothesis, there
would be,
for each common disease, a specific gene that is clearly responsible.
Such hopes were dashed over the coming years,
and, today, much research starts from the hypothesis that numerous genes
are differentially expressed in affected patients \cite{Goldstein}, but with individually
small effect sizes~\cite{Dai,CDRVH}. HC, with its emphasis on detecting rare
and weak effects,
seems well suited to this new environment.

\subsubsection{Two worked examples}
\label{subsec:realexample1}
Let's apply HC to the two gene microarray data sets. For each in turn,
let $x_{ij}$ denote the
expression level for the $i$th sample and the $j$th gene,
$1 \leq i \leq n$, $1 \leq j \leq p$. Let $C$ and $D$ be the set of
indices of samples from the
training set and the test set, respectively. For notational consistency
with later sections,
we only use the data in the training set, but using the whole data
gives similar results.

%
\begin{figure*}

\includegraphics{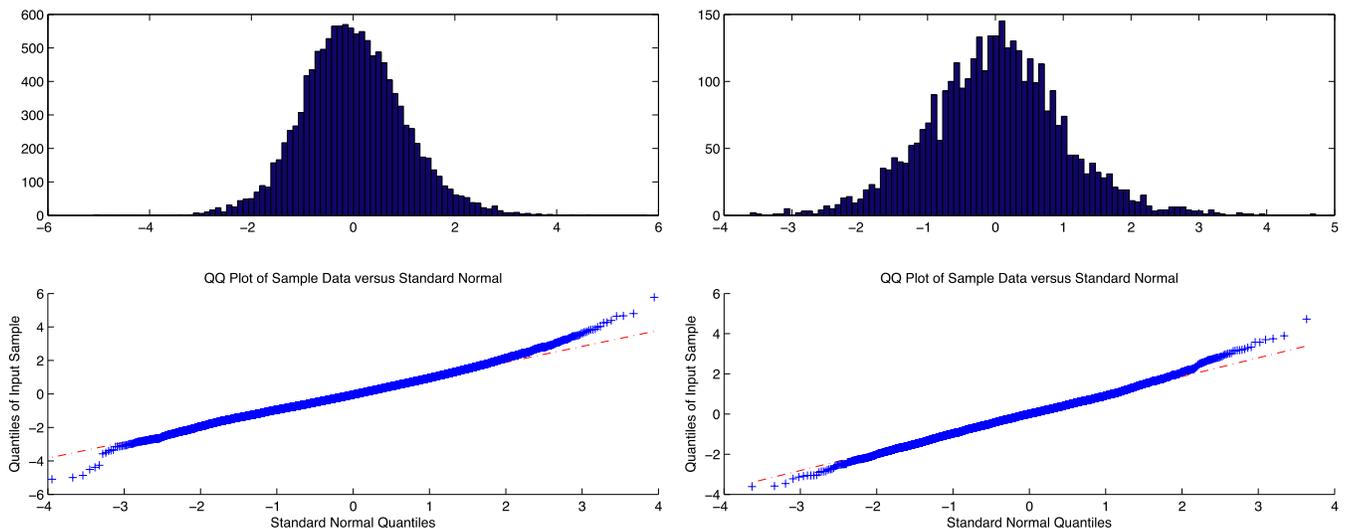}

\caption{Left column: histogram (top) and $qq$-plot (bottom) of
$Z = (Z_1, Z_2, \ldots, Z_p)'$ for the lung cancer data set.
Right: corresponding plots for the leukemia data set.}
\label{fig:cancersummary}
\end{figure*}

Write $C = C_1 \cup C_2$,
where $C_1$ and $C_2$ are the sets of indices
of the training samples from classes~$1$ and $2$, respectively.
Fix $1 \leq j \leq p$. Let $\bar{x}_{jk} = \frac{1}{|C_k|} \sum_{j
\in
C_k} x_{ij}$
denote the average expression value of gene $j$ for all samples in
class $k$, $k = 1, 2$,
and $s_j^2 = \frac{1}{(|C| -2)} [\sum_{i \in C_1} (x_{ij} - \bar
{x}_{j1})^2 + \sum_{i \in C_2} (x_{ij} - \bar{x}_{j2})^2]$ the pooled
variance.
Define the $t$-like statistic
\[
z_j^* = \frac{1}{\sqrt{1/|C_1| + 1/|C_2|}} \frac{\bar{x}_{j1} -
\bar
{x}_{j2}}{s_j},\quad 1 \leq j \leq p.
\]

In the null case, if the data $\{ x_{ij} \}_{1 \leq i \leq n, 1 \leq j
\leq p}$ are independent and identically distributed
across different genes, and, for each gene $j$, $\{x_{ij} \}_{i = 1}^n$
are normal samples with the same variance in each class, then $z_j^*$
has the Student's\vspace*{1pt} $t$-distribution with $df = |C| - 2$ when gene $j$ is
not differentially expressed.
Using this, we may calculate individual $P$-values for each gene and
apply HC.
However, as pointed out in Efron~\cite{EfronNull}, a problem one
frequently encounters in analyzing gene
microarray data is the so-called \textit{discrepancy between the empirical
null and theoretical null},
meaning that there is a gap between the aforementioned $t$-distribution
and the empirical null distribution
associated with $\{ z_j^* \}_{j = 1}^p$. This gap might be caused by
unsuspected between-gene variance components or other factors. We
follow Efron's suggestion and standardize $z_j^*$:
%
\begin{equation}
\label{renorm} Z_j = \frac{z_j^* - \bar{z}^*}{\mathrm{sd}(z^*)},\quad  1 \leq j \leq p,
\end{equation}
where $\bar{z}^*$ and $\mathrm{sd}(z^*)$ represent the empirical mean and
standard deviation associated with $\{z_j^*\}_{j = 1}^p$, respectively.
We call $Z_1, Z_2, \ldots, Z_p$
the standardized $Z$-scores, and believe that in the purely null case
they are approximately normally distributed.
In Figure~\ref{fig:cancersummary}, we show the histogram (top) and the
$qq$-plot (bottom) associated with the standardized
$Z$-scores. The figure suggests that, for both data sets, the
standardization in (\ref{renorm}) is effective.

We apply HC to $\{Z_j \}_{j = 1}^p$ for both the leukemia and the lung
cancer data, where the individual (two-sided) $P$-values are obtained
assuming $Z_j \sim N(0,1)$ if the $j$th gene is not differentially expressed.
The resulting HC scores are $6.1057$ and $13.3025$ in two cases.
The $P$-values associated with the scores (computed by numerical
simulations) are $\approx5 \times10^{-5}$ and $< 10^{-5}$, respectively.
They suggest the definite presence of signals, sparsely scattered in
the $Z$-vector.
And, indeed, the qqplots exhibit a visible ``curving away'' from the
identity lines.

An alternative approach to computing these two $P$-values uses random shuffles.
Denote the data matrix by $X = (x_{ij})_{1 \leq i \leq n, 1 \leq j \leq p}$.
First, we randomly shuffle the rows of $X$, independently between
different columns. For the permuted data, we follow the steps above and
calculate
the standardized $Z$-vector. Due to our shuffling, the signals wash out and
the $Z$-vector can be viewed as containing no effect
signals. Next, apply HC to the $Z$-vector,
obtaining an HC statistic
specific to that shuffle. Repeat the whole process for $1000$
independent shuffles.
As a result, we have $1001$ HC scores:
one, based on the original data matrix; all others, based on shuffles.
The shuffles are used to calculate the $P$-value of the original HC score.
For the (training) leukemia data and the (training) lung cancer data,
the resultant $P$-values are approximately $0.01$ and $ < 0.001$, respectively.
Both data sets have standardized $Z$-vectors exhibiting very subtle
departures from the
null hypothesis, but still permit reliable rejection of the null hypothesis.


\subsubsection{Applications to genome-wide association stu\-dy}
\label{subsec:GWAS1}

By now the literature of Genome-wide association studies (GWAS) has several
publications applying HC or its relatives.
Parkhomenko {et al}. \cite{GWAS1} used HC to detect modest
genetic effects
in a genome-wide study of rheumatoid arthritis. It was further
suggested in Martin {et al}. \cite{Martin}
that the implementation of HC in GWAS may provide evidence for the
presence of additional remaining SNPs modestly associated with the trait.

Sabatti {et al}. \cite{Sabatti} used HC in a GWAS for
metabolic traits.
HC enabled the authors to quantify the strength of the overall genetic
signal for each of the nine traits (Triglycerides, HDL$, \ldots$)
they were interested in, where, to deal with the possible dependence
caused by Linkage Disequilibrium (LD) between SNPs,
they computed individual $P$-values by permutations. See also De la
Cruz {et al}. \cite{delaCruz} where the authors
considered the problem of testing whether there are associated markers
in a given region or a given set of markers,
with applications to analysis of a SNP data set on Crohn's disease.

Wu {et al.} \cite{Wu} adapted HC for detecting rare and weak
genetic signals using
the information of LD. He and Wu \cite{GWAS2} used HC and innovated HC
for signal detection for large-scale exonic single-nucleotide
polymorphism data, and suggested modifications of HC in such settings.

Motivated by GWAS, Mukherjee {et al}. \cite{Lin} considered the
signal detection problem
using logistic regression coefficients rather than 2-sample $Z$-scores,
and discovered an interesting relationship between the sample size and the
detectability when both response variable and design variables are discrete.
See Section~\ref{subsec:corr} for more discussion on signal detection
problems associated
with regression models. To address applications in GWAS, Roeder and
Wasserman \cite{RW}
made an interesting connection between HC and weighted hypothesis testing.

\subsubsection{Applications to DNA copy number variation}
Computational biology continues to innovate and GWAS is no longer the
only game in town.
DNA Copy Number Variation (CNV)
data grew rapidly in importance after the GWAS era began,
and today provide an important window on genetic structural variation.
Jeng et al. \cite{Jeng1,Jeng2} applied HC-style thinking to CNV data
and proposed a new method called
\textit{Proportion Adaptive Segment Selection (PASS)}.
PASS can be viewed as a two-way screening procedure for genomic data,
which targets both the signal sparsity across different features (SNPs)
and the sparsity
across different subjects---so-called rare variation in genomics.

\subsection{Applications to Cosmology and Astronomy}
HC has been applied in several modern experiments in Astronomy and
Cosmology, where typically the
experiment produces data which can be interpreted as images (of a kind)
and where there
is a well-defined null hypothesis, whose overthrow would be considered
a shattering event.

Studies of the Cosmic Microwave Background (CMB) offer several examples.
CMB is a relic of radiation emitted when the Universe was about 370,000
years old.
In the simplest ``inflation'' models,
CMB temperature fluctuations should behave as a realization
of a zero-mean Gaussian random variable in each pixel.
The resulting Gaussian field (on the sphere) is completely determined by
its power spectrum. In recent decades, a large number of studies have
been devoted to the subject of detecting non-Gaussian signatures
(hot spots, cold spots, excess kurtosis$, \ldots$) in the CMB.

Jin {et al.} \cite{Starck}, and Cayon {et al.} \cite{Cayon1}
(see also \cite{Cruz,Cayon2}), applied
HC to standardized wavelet coefficients of CMB data
from the Wilkinson Microwave Anisotropy Probe (WMAP).
HC would be sensitive to a small collection of such coefficients
departing from the standard null, without requiring
that individual coefficients depart in a pronounced way.
Compared to the kurtosis-based non-Gaussianity detector (widely used
in cosmology when the departure
from Gaussianity is in the not-very-extreme tails), HC
showed superior power and sensitivity, and pointed, in particular,
to the \textit{cold spot} centered at galactic coordinate (longitude,
latitude) = $(207.8^{\circ}, -56.3^{\circ})$.
In \cite{Vielva}, Vielva reviews the cold spot detection problem and
shows that HC rejects
Gaussianity, confirming earlier detections by other methods.

Gravitational weak lensing calculations measure the
distortion of background galaxies supposedly caused by intervening large-scale
structure. Pires {et al.} \cite{Pires} applied many non-Gaussianity
detectors to weak lensing data,
including the empirical Skewness, the empirical Kurtosis and HC, and
showed that
HC is competitive, while of course being more specifically focused on
excess of observations in the tails
of the distribution.

Most recently, Bennett {et al}. \cite{Bennett} applied the HC ideas
to the
problem of Gravitational Wave detection.
They use HC as a second-pass method operating on $F$-statistic and
$C$-statistics (see \cite{Bennett} for details) for a monochromatic
periodic source in a binary system;
such statistics contain a large number of relatively weak signals
spread irregularly
across many frequency bands.
They use a modified form of HC, which is both sensitive and robust, and
offer a noticeable increase in the detection power (e.g., a $30\%$
increase in detectability for a phase-wandering source over multiple
time intervals).


\subsection{Applications to Disease Surveillance and Local Anomaly Detection}
In disease surveillance, we have aggregated count data $c_i$
representing cases of
a certain disease (e.g., influenza) by the $i$th spatial region (e.g.,
zip code), $1 \leq i \leq N$.
When disease breaks out, the counts will have elevated values in one or
a few
small geographical regions.
Neill and Lingwall \cite{Neill0,Neill} use HC for disease surveillance and
spatio-temporal cluster detection: they suppose we have historical
counts for each spatial location measured over time $t = 1, 2, \ldots, T$.
The $P$-value of $c_i$ is calculated by $(d_i + 1)/(T+1)$,
where $d_i$ is the number of historical counts larger than $c_i$ at the
$i$th location.

Disease outbreak detection is a special case of \textit{local anomaly}
detection
as studied in Saligrama and Zhao~\cite{anomaly}.
Suppose we have a graph $G = (V, E)$ with usual graph metric, where a
random variable
is associated with each node.
A simple scenario of \textit{local anomaly} they consider assumes
that, for all nodes outside the anomaly, the associated random
variables have the same density $f_0$,
and, for nodes inside the anomaly, the associated density is different
from $f_0$.
Saligrama and Zhao \cite{anomaly} investigate several models and
statistics for local anomaly detection; HC is found to be competitive
in this setting.

\subsection{Estimating the Proportion of Non-Null Effects}

As presented so far, HC offers a test statistic.
In the setting of Section~\ref{sec:AoS},
we sample $X_i$ from the two-component mixture
%
\begin{eqnarray}
\label{estepsmodel}\quad X_i \stackrel{\mathrm{i.i.d.}} {\sim} (1 - \eps)
N(0,1) + \eps N(\tau, 1),
\nonumber
\\[-8pt]
\\[-8pt]
 \eqntext{1 \leq i \leq N.}
\end{eqnarray}
The detection problem which HC addresses
involves testing $H_0^{(N)}\dvtx \eps= 0$ versus $H_1^{(N)}\dvtx \eps> 0$.
Alternatively, one could estimate the mixing proportion $\eps$.
Motivated by a study of Kuiper Belt Objects (KBO) (e.g., \cite{Rice}),
Cai {et al}. \cite{CJL} (see also \cite{Rice}) developed HC
into an
estimator for $\eps$, focusing on the regime where $\eps> 0$ is very small.\footnote{Among the many competing methods,
we mention just \cite{Wellner1}.
Let $\pi_i = 1 - \Phi(X_i)$ be as in model (\ref{estepsmodel});
then $\pi_i$ are i.i.d. samples from the density
$f_{\eps, \tau}(x) = (1 - \eps) + \eps g_{\eps, \tau}(x)$, where
$g_{\eps, \tau}(x)$ is monotone decreasing in $0 < x < 1$ and is
unbounded at $0$.
Balabdaoui {et al.} \cite{Wellner1} studied the behavior of the
maximum likelihood estimator of $f_{\eps, \tau}(x)$ at $0$,
and used it to derive an alternative estimator for $\eps$. }

In the growing literature of large-scale multiple testing, the problem of
estimating the proportion of nonnull effects has attracted
considerable attention in the past decade,
though sometimes with different goals. For example, rather than knowing
the true proportion
of nonzero effects, one might only want to estimate the largest proportion
within which the false discovery rate can be controlled.
The literature along this line connects to the work of Benjamini and
Hochberg \cite{BH95} on
controlling False Discovery Rate (FDR), and Efron \cite{EfronNull} on
controlling the local FDR in gene microarray studies.
See \cite{CaiJin2010,JinCai2007,Jin2008} and references therein.

\subsection{Statistics with HC-Like Constructions}
\label{sec:HCLike}

HC can be viewed as a measure of the goodness of fit between two
distributions, namely, between the distribution $F_N$ of the empirical
$P$-values
and the model uniform distribution $F_0$.
In this viewpoint, HC is effectively computing the distance measure
%
\begin{eqnarray}
\label{HCGF1} &&\mu_1( F_N, F_0)
\nonumber
\\[-8pt]
\\[-8pt]
\nonumber
&&\quad= \sqrt{N}
\cdot\max_{i=1}^{N} \frac{|F_N(i/N) -
F_0(i/N)|}{\sqrt{F_0(i/N)(1-F_0(i/N))}}
\end{eqnarray}
(or, more properly, a restricted form, where $i=1$ and $i > \alpha_0 N$
are omitted in the maximum) or the reverse
%
\begin{eqnarray}
\label{HCGF2} &&\mu_2( F_N , F_0)
\nonumber
\\[-8pt]
\\[-8pt]
\nonumber
&&\quad = \sqrt{N}
\cdot\max_{i=1}^{N} \frac{|i/N -
F_0(\pi
_{(i)})|}{\sqrt{F_0(\pi_{(i)}) (1 - F_0(\pi_{(i)}))}}.
\end{eqnarray}
We call these the \textit{theoretically standardized} and
\textit{empirically standardized} goodness of fit, respectively.

To understand HC, then, one might consider
how it differs from other measures of the discrepancy
between two distributions. HC includes
the element of standardization, which for many readers will
suggest comparison to the Anderson--Darling statistic \cite{AndersonDarling}:
\[
A( F_N , F_0) = N \cdot\int\frac{|F_N(x) -
F_0(x)|^2}{F_0(x)(1-F_0(x))} \,dx.
\]
HC, however, involves maximization rather than integration,
which makes it a kind of weighted Kolmogorov--Smirnov statistic.
Jager and Wellner \cite{Wellner2004} investigated the limiting
distribution of a class of
weighted Kolmogorov statistics, including HC as a special case.

Another perspective is to view the $P$-values underlying HC
as obtained from the normal approximation to a one-sample test for a
known binomial
proportion, and to consider instead the exact test based on likelihood
ratios, or asymptotic tests
based instead on KL divergence between the binomial with parameter $\pi
_i$ and the binomial with
parameter $i/N$. This perspective reveals a similarity of HC to the
\textit{Berk--Jones} (BJ)
statistic \cite{BerkJones}. The similarity was carefully studied in~\cite{DJ04},
Section~1.6; see details therein. Using the divergence
$D(p_0,p_1) = p_0 \log(p_0/p_1) + (1-p_0) \log((1-p_0)/(1-p_1))$,
the Berk--Jones statistic can be written as
\[
\mathrm{BJ} = \max_{i=1}^N N \cdot D(\pi_i,i/N).
\]
Wellner and Koltchinskii \cite{WellnerKoltchinskii}
derive the limiting distribution of the Berk--Jones statistic, finding
that it
shares many theoretical properties in common with HC.

In \cite{Wellner2007}, Jager and Wellner introduced a new family of
goodness-of-fit tests based on the $\phi$-divergence, including HC as a
special case,
and showed all such tests achieve the optimal detection boundary in
\cite{DJ04} (see the discussion below in Section~\ref{sec:phase}).

Reintroducing the element of integration found in the
Anderson--Darling statistic,
Walther \cite{Walther} proposed an \textit{Average Likelihood Ratio}
(ALR) approach.
If $\mathrm{LR}_{i,N} $ denotes the usual likelihood
ratio for a one-sided test of the binomial proportion, ALR takes the form
\begin{eqnarray*}
\mathrm{ALR} &= &\sum_{i=1}^{\alpha_0 N} w_{i,N}
\mathrm{LR}_{i,N}; \\
\mathrm{LR}_{i,N} &\equiv& \exp\bigl( N \max\bigl\{ D(
\pi_i,i/N),0 \bigr\}\bigr),
\end{eqnarray*}
with weights $w_{i,N} = (2 i \log(N/3))^{-1}$.
Walther shows that ALR compares favorably with
HC and BJ for finite sample performance, while having similar
asymptotic properties under the ARW model discussed below.
See \cite{DJ04,Tukey76} for more discussions on the relative merits of HC, BJ
and ALR.

Additionally, as a measure of goodness of fit, HC is closely
related to other goodness-of-fit tests, motivated, however, by the goal
of optimal detection of
presence of mixture components representing rare/weak signals.
We remark that the pontogram of Kendall and Kendall \cite{Kendall}
is an instance of HC, applied to a special set of $P$-values.

Gontscharuk {et al.} \cite{Gontscharuk}
introduced the notion of \textit{local levels} for goodness-of-fit tests
and studied the
asymptotic behavior when applying the framework to one version of HC;
for HC, the
local level associated with $\mathrm{HC}_{N, i}$ and a critical value $\chi$
roughly translates to
$P[\mathrm{HC}_{N, i} \geq\chi]$.


\subsection{Connection to FDR-Controlling Methods}
\label{subsec:FDR}
HC is connected to Benjamini and Hochberg's (BH) False Discovery Rate
(FDR) control method in large-scale multiple testing \cite{BH95}.
Given $N$ uncorrelated tests where $\pi_{(1)} < \pi_{(2)} < \cdots<
\pi_{(N)}$
are the sorted $P$-values, introduce the ratios
\[
r_k = \pi_{(k)} / (k/N),\quad 1 \leq k \leq N.
\]
Given a prescribed level $0 < q < 1$ (e.g., $q = 5\%$), let
$k = k_q^{\mathrm{FDR}}$ be the largest index such that $r_k \leq q$.
BH's procedure rejects all tests whose $P$-values are among
the $k_q^{\mathrm{FDR}}$ smallest, and accepts all others. The procedure
controls the FDR in that
the expected fraction of false discoveries is no greater than $q$.

The contrast between FDR control method and HC can be captured in a few
simple slogans.
We think of the BH procedure as targeting \textit{rare but strong} signals,
with the main goal to select the few strong signals embedded in a long
list of null signals,
without making too many false selections. HC targets the more delicate
regime where the signals are \textit{rare and weak}.
In the rare/weak setting, the signals and the noise may be almost
indistinguishable; and while
the BH procedure still controls the FDR, it yields very few
discoveries. In this case,
a more reasonable goal is to test whether any signals exist without
demanding that we properly identify them all; this is what HC is specifically
designed for.
See also Benjamini \cite{FDR1}.

HC is also intimately connected to the problem of constructing
confidence bands for the
\textit{False Discovery Proportion} (FDP). See Cai et al. \cite{CJL},
Ge and Li \cite{Ge}, and de Una-Alvarez \cite{Una}.

\subsection{Innovated HC for Detecting Sparse Mixtures in Colored Noise}
So far, the underlying $P$-values were always assumed independent.
Dai {et al.} \cite{Dai} pointed out the
importance of the correlated case for genetics and genomics;
we suppose many other application areas have similar concerns.
Hall and Jin \cite{HJ}
showed that directly using HC in such cases
could be unsatisfactory, especially under strong correlations. Hall and
Jin \cite{InnovatedHC}
pointed out that correlations (when known or accurately estimated)
need not be a nuisance or curse,
but could sometimes be a blessing if used properly.
They proposed \textit{Innovated Higher Criticism},
which applies HC in a transformed coordinate system;
in analogy to time series theory, Hall and Jin called this
the innovations domain.
Innovated HC was shown to be successful when the correlation matrix
associated with the noise entries has polynomial off-diagonal decay.

\subsection{Signal Detection Problem Associated with Regression Models}
\label{subsec:corr}
Suppose we observe an $n \times1$ vector $Y$ which satisfies
a linear regression model
\[
Y = X \beta+ z, \quad z \sim N(0, I_n),
\]
where $X$ is the $n \times N$ design matrix, $\beta$ is the $N \times
1$ vector of regression coefficients, and $z$ is the noise vector.
The problem of interest is now to test whether all regression
coefficients $\beta_i$ are $0$ or a small fraction of
them is nonzero. The setting considered in \cite{HJ,InnovatedHC} is a
special case, where
the number of variables $N$ is the same as the sample size $n$.

Arias-Castro {et al.} \cite{Castro2} and Ingster, Tsybakov and
Verzelen \cite{IPV} considered the more general case where $N$ is much
larger than $n$.
The main message is that, under some conditions, what has been previously
established for the Gaussian sequence model extends to high-dimensional
linear regression.
Motivated by GWAS, Mukherjee {et al}. \cite{Lin} considered a
similar problem
with binary response logistic regression.
They exposed interesting new phenomena governing the detectability of
nonnull $\beta$
when both response variable and design variables are discrete.

Meinshausen \cite{Meinshausen} considers the problem of
variable selection associated with a linear model.
Adapting HC to the case of correlated noise with unknown variance,
he uses the resultant method for hierarchical testing of variable importance.
Charbonnier \cite{Charbonnier} generalizes HC from a one-sample testing
problem to a two-sample testing problem. It considers two linear models
and tries to
test if the regression coefficient vectors are the same. Also related is
Suleiman and Ferrari \cite{Suleiman}, where the authors
use constrained likelihood ratios for detecting sparse signals in
highly noisy 3D data.

\subsection{Signal Detection when Noise Distribution is Unknown/Non-Gaussian}
In models (\ref{mixture1})--(\ref{mixture2}), the noise entries are i.i.d.
samples from $N(0,1)$.
In many applications, the noise distribution is unknown and is probably
non-Gaussian.
To use HC for such settings, we need an approach to computing
$P$-values $\pi_i$, $1 \leq i \leq N$.

Delaigle and Hall \cite{Delaigle} and Delaigle {et al}. \cite
{Robustness}
addressed this problem in the settings where the data are arranged in a
$2$-D array
$\{X(i,j)\}$, $1 \leq i \leq n$, $1 \leq j \leq N$.
In this array, different columns are independent, and entries
in the $j$th column are i.i.d. samples from a distribution $F_j$
which is unknown and presumably non-Gaussian.
We need to associate a $P$-value with each column.
A conventional approach is to compute the $P$-value using the
Student's $t$-statistic. However, when $F_j$ are non-Gaussian, the
$P$-values may not be accurate enough,
and the authors propose to correct the $P$-values with bootstrapping.
A similar setting is considered by Greenshtein and Park \cite
{Greenshtein} and by Liu and Shao \cite{ShaoQM}. The first paper
proposes a modified
Anderson--Darling statistic and shows that, in certain settings, the
proposed approach
may have advantages over HC in the presence of non-Gaussianity.
The second paper proposes a test based on extreme values of
Hotelling's $T^2$, and studies the case where the sparse signals appear
in groups and the underlying distributions are not necessarily normal.

\subsection{Detecting Sparse Mixtures more Generally}
More generally, the problem of detecting sparse mixtures considers hypotheses
\begin{eqnarray*}
&&H_0^{(N)}\dvtx\quad X_i \stackrel{\mathrm{i.i.d.}}
{\sim} F,\quad vs.\\
&& H_1^{(N)}\dvtx\quad X_i \stackrel{
\mathrm{i.i.d.}} {\sim} (1 - \eps) F + \eps G,
\end{eqnarray*}
where $\eps\in(0,1)$ is small and $F$ and $G$ are two distributions
that are presumably different; $(\eps, F, G)$ may depend on $N$. In
Donoho and Jin \cite{DJ04}, $F = N(0,1)$ and $G = N(\tau, 1)$ for some
$\tau> 0$.

Cai {et al}. \cite{CJJ} considered the case where $F = N(0,1)$ and
$G = N(\tau, \sigma^2)$, so the mixture in the alternative hypothesis
is not only heterogeneous but also heteroscedastic, and $\sigma$ models
the heteroscedasticity. They found that $\sigma$ has a surprising
phase-change effect over the detection problem. The heteroscedastic
model is also considered in Bogdan {et al}. \cite{Bogdan1} and
Bogdan {et al}. \cite{Bogdan2} from a
Bayesian perspective.
Park and Ghosh \cite{Ghosh} gave a nice review on recent topics on
multiple testing where HC is discussed in detail.

Cai and Wu \cite{CaiWu} extend the study to the more general case where
$F = N(0,1)$ and $G$ is a Gaussian location mixture with a general
mixing distribution, and study the detection boundary as well as the
detectability of HC.

Arias-Castro and Wang \cite{WangMeng} investigate the case where $F$ is
\textit{unknown} but symmetric,
and develop distribution-free tests to tackle several interesting
problems, including
that of testing of symmetry.

In addition, Gayraud and Ingster \cite{Ingstervariable} consider the
problem of detecting sparse mixtures in the
functional setting, and show that the HC statistic continues to be
successful in the very sparse case.
Laurent {et al}. \cite{Laurent} considered the problem of testing
whether the samples $X_i$ come
from a single normal or a mixture of two normals with different means
(both means are unknown).

In a closely related setting, Addario-Berry {et al}. \cite{Lugosi}
and Arias-Castro {et al}. \cite{Erycluster}
considered structured signals, forming clusters in geometric shapes
that are unknown to us.
The setting is closely related to the one considered in \cite{InnovatedHC},
Section~6.
Haupt et al. \cite{Haupt1,Haupt2} considered a more complicated setting
where an adaptive sample scheme is available, where we can do inference
and collect
data in an alternating order.

\section{Higher Criticism for Feature Selection} \label{sec:class}

Higher Criticism has applications far beyond the testing of a global
null hypothesis.

Consider a classification problem where we have training samples
$(X_i, Y_i)$, $1 \leq i \leq n$, from two different classes.
We denote $X_i$ by the feature vectors and $Y_i = \pm1$ the class
labels. For simplicity,
we assume two classes are equally likely and the feature vectors $X_i
\in R^p$ are Gaussian distributed with identical covariances, so that,
after a standardizing transformation,
the feature vector $X_i \sim N(Y_i \cdot\mu, I_p)$, with vector $\mu
$ being the contrast mean and $I_p$ the $p \times p$ identity matrix.
Given a fresh feature vector $X$, the goal is to predict the
associated class label $Y \in\{-1, 1\}$.

We are primarily interested in the case where $p \gg n$ and where the
contrast mean vector $\mu$ is unknown
but has nonzero coordinates that are both rare and weak. That is, only
a small fraction of coordinates of $\mu$ is nonzero, and
each nonzero coordinate is individually small and contributes weakly to
the classification decision.

In the classical $p < n$ setting, consider traditional Fisher linear
discriminant analysis (LDA).
Letting $w = (w(j), 1 \leq j \leq p)$ denote a sequence of feature
weights, Fisher's LDA takes the form
\[
L(X) = \sum_{j =1}^p w(j) X(j).
\]
It is well known that the optimal weight vector $w \propto\mu$, but
unfortunately
$\mu$ is unknown to us and in the $p > n$ case can be
hard to estimate, especially when the nonzero coordinates of $\mu$ are
rare and weak; in that case, the empirical estimate
$\bar{X}$ is noisy in every coordinate, and only a few coordinates
``stick out'' from the noise background.

Feature selection (i.e., selecting a small fraction of the available
features for classification) is a standard approach to attack the
challenges above.
Define a vector of feature scores
%
\begin{equation}
\label{Zvector} Z = \frac{1}{\sqrt{n}} \sum_{i = 1}^n
(Y_i \cdot X_i);
\end{equation}
this contains the evidence in favor of each feature's significance.
We will select a subgroup of features for our classifier,
using hard thresholding of the feature scores.
For a threshold value $t > 0$ still to be determined,
define the hard threshold function
\[
w_t(z) = \mathrm{sgn}(z) \cdot1_{\{|z| > t\}},
\]
which selects the features having sufficiently large evidence
and preserves the sign of such feature scores.
The post-feature-selection Fisher's LDA rule is then
\[
L_t(X) = \sum_{j = 1}^p
w_t(j) X(j),
\]
and we simply classify $Y$ as $\pm1$ according to $L_t(X) \gtrless0$.
This is related to the modified HC in \cite{Chen}, but there the focus
is on signal detection instead of feature selection.


How should we set the threshold $t$? Consider \textit{HC feature selection},
where a simple variant of HC is used to set the threshold.
To apply HC to feature selection,
we fix $\alpha_0 \in(0, 1/2]$ and follow three steps (to be consistent
with OHC described in Section~\ref{subsec:HCbasic}, we switch back from
$p$ to $N$;
note that $N = p$ in this section):
\begin{itemize}
\item Calculate a (two-sided) $P$-value
$\pi_j = P\{ |N(0,1)| \geq|Z(j)|\}$ for each $1 \leq j \leq N$.
\item Sort the $P$-values into ascending order: $\pi_{(1)} < \pi_{(2)}
< \cdots< \pi_{(N)}$.
\item Define the \textit{Higher Criticism feature scores} by
%
\begin{eqnarray}
&&\mathrm{HC}(i; \pi_{(i)}) = \sqrt{N} \frac{i/N - \pi_{(i)}}{\sqrt{(i/N)(1 -
i/N)}},
\nonumber
\\[-8pt]
\\[-8pt]
\eqntext{1 \leq i \leq N.}
\end{eqnarray}
Obtain the maximizing index of $\mathrm{HC}(i; \pi_{(i)})$:
\[
\hat{i}^{\mathrm{HC}} = \margmax_{\{1 \leq i \leq\alpha_0 \cdot N \}} \bigl\{\mathrm{HC}(i,
\pi_{(i)}) \bigr\}.
\]
The \textit{Higher Criticism threshold (HCT)} for feature selection is then
\[
\that_N^{\mathrm{HC}} = \that_N^{\mathrm{HC}}(Z_1,
Z_2, \ldots, Z_N; \alpha_0, n) =
|Z|_{\hat{i}^{\mathrm{HC}}}.
\]
\end{itemize}
In modern high-throughput settings where {a priori} relatively
few features are likely to be useful,
we set $\alpha_0 = 0.10$.\footnote{In practice, HCT is relatively insensitive to different
choices of~$\alpha_0$.}$^{,}$\footnote{Note the denominator of the HC objective function is
different from the denominator used earlier, in testing, although the
spirit is similar.
The difference is analogous to the one between the two goodness-of-fit tests
(\ref{HCGF1}) and (\ref{HCGF2}).}
See \cite{DJ09} for explanation.

Once the threshold is decided, LDA with HC feature selection is
\begin{eqnarray}
L_{\mathrm{HC}}(X) = \sum_{j = 1}^p
w_{\mathrm{HC}}(j) X(j),
\nonumber
\\
\eqntext{\mbox{where } \displaystyle w_{\mathrm{HC}}(j) = \mbox{sgn}\bigl(Z(j)
\bigr) 1\bigl\{\bigl |Z(j)\bigr| \geq\that_p^{\mathrm{HC}} \bigr\},}
\end{eqnarray}
and the HCT trained classification rule will classify $Y = \pm1$
according to $L_{\mathrm{HC}}(X) \gtrless0$.

The classifier above is a computationally inexpensive approach,
especially when compared to resampling-based methods (such as
cross-validations, boosting, etc.). This gives HC a lot of
computational advantage in the now very relevant ``Big Data'' settings.

%
\begin{figure*}

\includegraphics{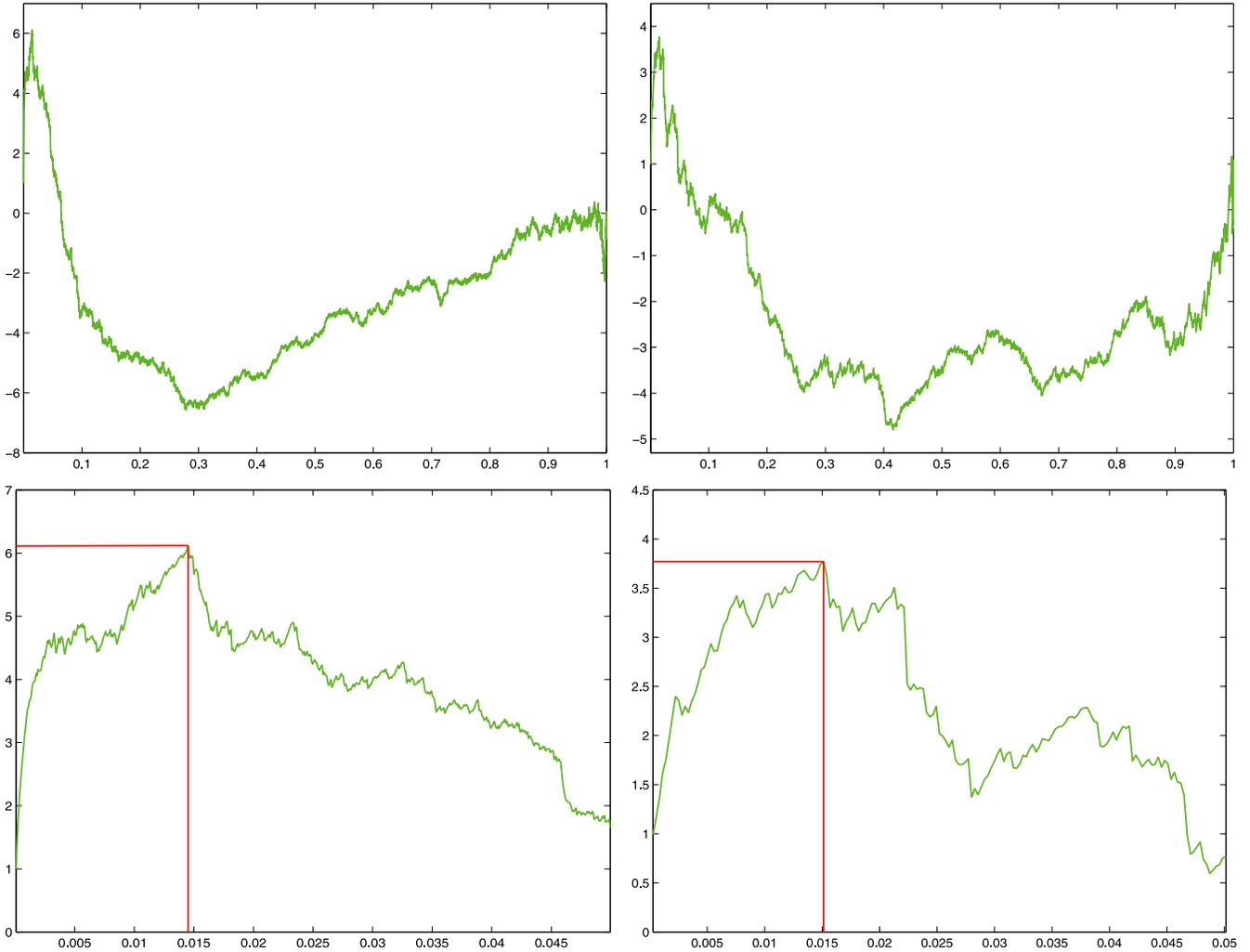}

\caption{Top row: plot of feature scores $\mathrm{HC}_{N,i}$ versus $i/N$ for
lung cancer data (left) and leukemia data (right). Bottom row:
enlargements of plots in the top row.}
\label{fig:cancerHC}
\end{figure*}

\subsection{Applications to Gene Microarray Data} \label{sec:application}
We now apply the HCT classification rule to the two microarray data
sets discussed earlier. Again, $Z_j$ is the standardized $Z$-score
associated with the $j$th gene, using all samples in the training set $C$.

%

\begin{figure*}

\includegraphics{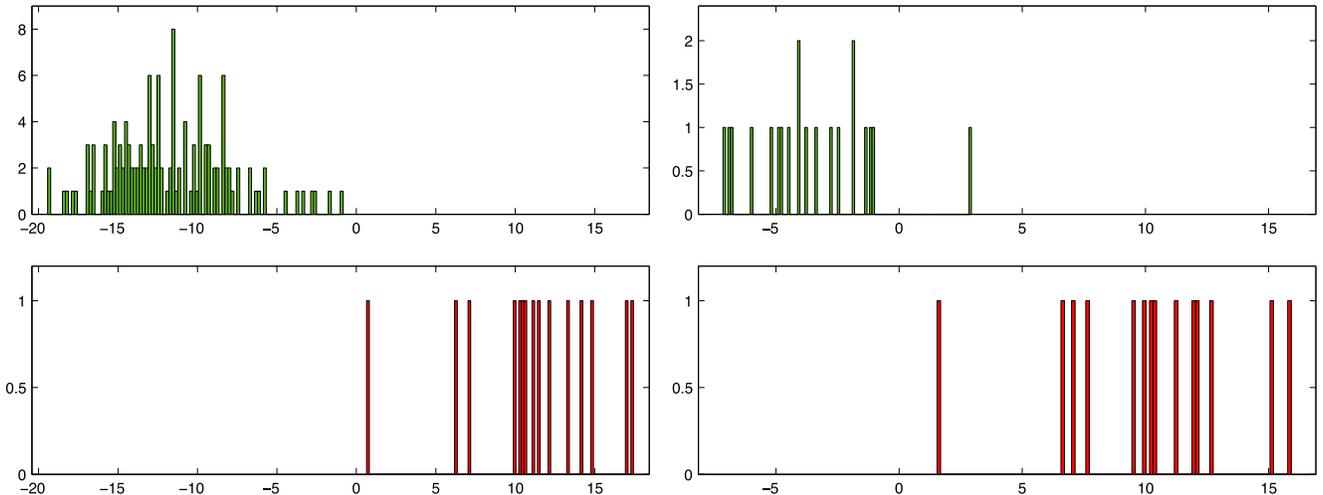}

\caption{Top row: histogram of the test scores corresponding to class
1. Bottom row: class 2. Left column: lung cancer data. Right column:
leukemia data.}
\label{fig:cancerscores}
\end{figure*}

First, we apply HC to $Z = (Z_1, Z_2, \ldots, Z_p)'$ to obtain the HC
threshold; this will also determine how many features we
keep for classification. The scores $\mathrm{HC}_{N, i}$ are displayed in
Figure~\ref{fig:cancerHC}.
For the lung cancer data,
the maximizing index is $\hat{i}^{\mathrm{HC}} = 182$, at which the HC score is $6.112$,
and we retain all $182$ genes with the largest $Z$-scores (in absolute
value) for classification
(equivalently, a gene is retained if and only if the $Z$-score exceeds
$\hat{t}_p^{\mathrm{HC}} = 2.65$).
For the leukemia data, $\hat{i}^{\mathrm{HC}} = 54$,
with HC score $3.771$ and the threshold $\hat{t}_p^{\mathrm{HC}} = 2.68$.

Next, for each sample $X_i$ in the test set $D$, we calculate
the HCT-based LDA score.
Recall that for any $i \in C$, the data associated with the $i$th
sample is $X_i = \{x_{ij}\}_{j = 1}^p$.
The HCT-based LDA score $\mathrm{lda}_{i} = \mathrm{lda}(X_i)$ is given by
\[
\mathrm{lda}_i = \sum_{j = 1}^p
w_{\mathrm{HC}}(j) \biggl( \frac{x_{ij} - \bar{x}_j}{s_j} \biggr),
\]
where we recall $w_{\mathrm{HC}}(j) = \mathrm{sgn}(Z_j) 1\{ |Z_j| \geq\that
_p^{\mathrm{HC}} \}$, and
$(\bar{x}_j, s_j, w_{\mathrm{HC}}(j))$ only depend on the training samples in $C$.
The scores $\{\mathrm{lda}_i\}_{i = 1}^n$ are displayed in
Figure~\ref{fig:cancerscores}, where we normalized each score by a
common factor of $1/ \sqrt{\hat{i}^{\mathrm{HC}}}$ for clarity.
The scores corresponding to class 1 are displayed
in the top row in green (ADCA for lung cancer data and ALL for leukemia),
and the scores for class 2 are displayed in the bottom row in red (the
left column displays lung cancer data and the right column displays
leukemia data).
For lung cancer data, LDA--HCT correctly classifies each sample.
For leukemia data, LDA--HCT correctly classifies each sample, with one exception:
sample $29$ in the test set
(number $67$ in the whole data set).

We employed these two data sets because they gave such a clear
illustration. In our previous paper \cite{DJ08}, we considered each
data in Dettling's well-known compendium, which includes the colon
cancer data and the prostate data.
The results were largely as good as or better than other classifiers,
many involving much fancier-sounding underlying principles.

\subsection{Threshold Choice by HC: Applications to Biomarker Selection}
\label{subsec:GWAS2}

Wehrens and Fannceschi \cite{Wehrens} used HC thresholding for
biomarker selection to analyze
metabolomics data from spiked apples.
They considered $P$-values that are calculated from Principal Component
scores and reported a marked
improvement in biomarker selection, compared to the standard selection
obtained by existing practices.
The paper concludes that HC thresholds can differ considerably from
current practice, so it is
no longer possible to blindly apply the selection thresholds
used historically; the data-specific cutoff values provided by HC open
the way to objective comparisons
between biomarker selection methods, not biased by arbitrary habitual
threshold choices.

\subsection{Comparison to Other Classification Approaches}
The LDA--HCT classifier is closely related to
other threshold-based LDA feature selection rules:
PAM by Tibshirani {et al}. \cite{Tibs} and FAIR by Fan and Fan
\cite
{FanFan}.
HCT picks the threshold based on feature $Z$-scores by
Higher Criticism, while the other methods
set this threshold differently. For the same data sets we
discussed earlier, the error rates for PAM and FAIR were reported in
\cite{FanFan,Tibs};
as it turns out, LDA--HCT has smaller error rates.

Comparisons with some of the more popular ``high-tech'' classifiers
(including Boosting \cite{Boosting},
SVM \cite{SVM} and Random Forests \cite{Breiman})
were reported in \cite{DJ08}. More complex methods usually need careful
tuning to perform well,
but HCT--LDA is very simple, both conceptually and computationally. When
used on the ensemble of standard
data sets published in Dettling, HCT--LDA happens to be minimax-regret
optimal: it suffers
the least performance loss, relative to the best method, across the ensemble.

Hall {et al}. \cite{Hall} apply HC for classification in a
different manner.
They view HC as a goodness-of-fit diagnostic. Their method first uses
the training vectors to obtain the empirical
distributions of each class, and then uses HC to tell which of these
distributions best fits each test vector.
They classify each test vector using the best-fitting class distribution.
While this rule is sensible, it turns out that
in a formal asymptotic analysis using the rare/weak model,
it is outperformed substantially by HCT--LDA.

\subsection{Connection to Feature Selection by Controlling Feature-FDR}
\label{subsec:featureFDR}
False Discovery Rate control methods offer a popular approach for
feature selection.
Fix $0 < q < 1$. $\mathrm{FDRT}_q$ selects features in a way so that
\[
\mbox{feature-FDR} \equiv E \biggl[\frac{\#\{\mbox{Falsely selected
features}\}}{\#\{\mbox{All selected features}\}} \biggr] \leq q.
\]
In the simple setting considered in Section~\ref{sec:class}, this can
be achieved by applying Benjamini--Hochberg's FDR controlling method to
all feature $P$-values.
The approach appeals to the common belief that, in order to have optimal
classification behavior, we should select features in a way so that the
feature-FDR stays small.

However, such beliefs have theoretical support only when signals are
rare/strong. In principle, the optimal $q$ associated with the
optimal classification behavior should depend on the underlying
distribution of the signals (e.g., sparsity and
signal strength); and when signals are rare/weak, the optimal FDR level
turns out to be much larger than $5\%$, and in some cases is close to~$1$.
In \cite{DJ09}, we studied the optimal level in an asymptotic rare/weak
setting and derived the leading asymptotics of
the optimal $\mathrm{FDR}$. In Section~\ref{subsec:phaseFDR} below we give more detail.

In several papers \cite{Strimmer,Klaus,Klaus2}, Strimmer and
collaborators compared the approach of feature selection by HCT
with both that of control of the FDR and that of control of the
False non-Discovery Rate (FNDR), analytically and also with synthetic
data and several real data sets on cancer gene microarray.
In their papers, they also compared the EBayes approach of Efron~\cite{Efron},
which presets an error rate threshold (say, $2.5\%$) and targets a
threshold where
the prediction error falls below the desired error rate.
Their numerical studies confirm the points explained above: HCT adapts
well to
different sparsity level and signal strengths, while the methods of
controlling FNDR and EBayes
do not perform as well (in the misclassification sense); and HC
typically selects
more false features than other approaches.
The goal of the HC feature selection, as we will see,
is to optimize the classification error, not to control the FDR. In fact,
\cite{Klaus}, Table~2, found that HCT had the best classification
performance for the cancer microarray data sets they investigated.

\subsection{Feature Selection by HCT When Features Are Correlated}
Above, we assumed the feature vector $X_i \sim N(Y_i \cdot\mu, I_p)$
for $Y_i = \pm1$.
A natural generalization is to assume $X_i \sim N(Y_i \cdot\mu,
\Sigma
)$, where $\Sigma= \Sigma_{p, p}$
is a unknown covariance matrix.
Two problems arise: how to estimate the precision matrix $\Omega=
\Sigma^{-1}$ and how to incorporate the estimated precision matrix into
the HCT classifier. In the latter, the key is to extend the idea of
threshold choice by HCT to the setting where not only the features are
correlated, but the covariance matrix is unknown and must be estimated.

The authors of \cite{HuangJin} address the first problem by proposing
\textit{Partial Correlation Screening (PCS)} as a new row-wise
approach to
estimating the precision matrix.
PCS starts by computing the $p \times p$ empirical scatter matrix $S =
(1/n) \sum_{i = 1}^n X_i X_i'$.
Assume the rows of $\Omega$ are relatively sparse.
To estimate a row of $\Omega$, the algorithm only needs to access
relatively few rows of~$S$.
For this reason, the method is able to cope with much larger $p$ (say,
$p = 10^4$)
than existing approaches (e.g., Bickel and Levina \cite{BLT}, glasso
\cite{glasso},
Neighborhood method \cite{MB} and CLIME \cite{CaiLuo}). \cite{FanHC}
addresses the second problem by combining the ideas in
Donoho and Jin \cite{DJ08} on threshold choices by HCT with those
in Hall and Jin \cite{HJ} on Innovated HC.
This combination injects
an estimate of $\Omega$ into the HCT classification method; it
is asymptotically optimal if $\Omega$ is sufficiently sparse and we
have a reasonably good estimate of $\Omega$ (e.g., \cite{CaiLuo}).

\section{Testing Problems About a Large Covariance Matrix}
In this section and the next, we briefly develop stylized applications
of HC to settings
which may seem initially far outside the original scope of the idea.
In each case, HC requires merely the
ability to compute a collection of $P$-values
for a collection of statistics
under an intersection null hypothesis.
This allows us to easily obtain HC-tests in
diverse settings.

Consider a data matrix $X = X_{n,p}$, where the rows of $X$ are i.i.d.
samples from
$N(0, \Sigma)$.
We are interested in testing $\Sigma= I_p$ versus the hypothesis that
$\Sigma$ contains a sub-structure.
First, we consider the case where the substructure is a small-size clique.
In Section~\ref{subsec:clique}, we approach the testing problem by
applying HC to the whole body of pairwise
empirical correlations and to the maximum row-wise correlation (for
each variable,
this is the maximum of each variable's correlations
with all other variables).
Second, in Section~\ref{sec:lowrank}, we consider the case where the matrix
$\Sigma= I + H$ follows the so-called \textit{spiked covariance model}
\cite{John02},
a low-rank perturbation of the identity matrix.
We apply HC to the eigenvalues
of the empirical covariance matrix.

\subsection{Detecting a Possible Clique in the Covariance Matrix}
\label{subsec:clique}
In this section, the global null hypothesis is $\Sigma= I_p$,
while the alternative is that $\Sigma$ contains a small clique.
Formally, $\Sigma$ can be written as $\Sigma= \Gamma\Sigma_0 \Gamma'$,
where $\Gamma$ is a permutation matrix, and, for an integer $1 \leq k
< p$ and $a \in[0, 1)$,
%
\begin{eqnarray}
\label{Definecliquesigma} &&\Sigma_0(i,j)
\nonumber
\\[-8pt]
\\[-8pt]
\nonumber
&&\quad = \cases{ %
1\{i = j \} + a \{ i \neq j \},&
$ \max\{i,j \} \leq k,$
\vspace*{2pt}\cr
1\{ i = j\}, &$\max\{i,j\} > k.$}
\end{eqnarray}
The parameter $a$ can take negative values as long as $\Sigma_0$
remains positive definite.

We suggest two different approaches for detecting the cliques using HC.
In each of the two approaches, the key is to obtain $P$-values.

In the first approach, we obtain individual $P$-values from pairwise
correlations.
In detail, write the data matrix $X = X_{n,p}$ as
\[
X = [x_1, x_2, \ldots, x_p].
\]
The pairwise correlation between the $i$th and $j$th variable is
\[
\rho_{ij} = \frac{(x_i, x_j)}{\|x_i\| \| x_j \| }.
\]
Recall that $t_{k}(0)$ denotes the central Student's $t$-distribution
with $df = k$.
The following lemma summarizes some basic properties of $\rho_{ij}$
\cite{Ruben}.
%
\begin{lemma} \label{lemma:pairwisecorr}
Suppose $\Sigma= I_p$. If $i \neq j$, then for any $\rho\in(-1, 1)$,
$P( \rho_{ij} \geq\rho) = P(t_{n-1}(0) \geq\sqrt{n-1} \rho/\break \sqrt{1
- \rho^2})$, Also,
if $(i, j) \neq(k, \ell)$, then $\rho_{ij}$ and $\rho_{k \ell}$ are
independent.
\end{lemma}
This says that the collection of random variables
\[
\{\rho_{ij}\dvtx 1 \leq i \leq j \leq p\}
\]
are pairwise independent (but not jointly independent). It can be
further shown that the correlation matrix between different $\rho_{ij}$
is very sparse,
so a simple but reasonable approach is to apply OHC to $\{\rho_{ij}\dvtx 1
\leq i \leq j \leq p\}$ directly. Numerically, the correlation between
different $\rho_{ij}$ will not significantly affect the performance of
OHC. On the other hand, since the correlation matrix
between $\rho_{ij}$ can be calculated explicitly,
a slightly more complicated method is to incorporate the correlation
structures into HC, following the idea of Innovated HC \cite{InnovatedHC}.

In the second approach, we obtain $P$-values from the maximum
correlation in each row:
\[
\rho_i^* = \max_{j \neq i} \rho_{ij}, \quad 1 \leq
i \leq p.
\]
%
%
\begin{lemma} \label{lemma:maxrho}
Suppose $\Sigma= I_p$. For $1 \leq i \leq p$ and $\rho\in(-1, 1)$,
\begin{eqnarray*}
P\bigl( \rho_i^* \leq\rho\bigr) &=& \biggl[P \biggl(t_{n-1}(0)
\leq\frac
{\sqrt
{(n-1)} \rho}{\sqrt{1 - \rho^2}} \biggr) \biggr]^{p-1} \\
&\equiv& F_{p,n}(
\rho).
\end{eqnarray*}
\end{lemma}
For a proof, see \cite{Ruben}, for example.
Let $\pi_i^* = F_{p,n} (\rho_i^*)$, so under the global null,
$\pi_i^* \sim\mathrm{Unif}(0,1)$.
We simply use these $P$-values in the standard HC framework.\footnote{$\pi_i^*$ are equi-correlated: for any $1 \leq i \neq j
\neq p$,
$\mathrm{Cov}(\pi_i^*, \pi_j^*) = c_0(n,p)$ for a small constant
$c_0(n,p)$ that does not depend on $i$ or $j$ and can be calculated numerically.
It can be shown that $c_0(p,n) = O(1/p)$, and the equi-correlation does
not have a major influence asymptotically.
Numerical study confirms that correcting for the equi-correlation only
has a negligible difference, so we only report results without the
correction. }

We conducted a small-scale simulation as follows. Fix $(p, n) = (1000,500)$.
We consider $5$ different combinations of $(k, a)$: $\{(1,0), (5,
0.25), (15, 0.2),\break  (45, 0.1), (135, 0.05)\}$. For each combination,
define $\Sigma_0$ as in (\ref{Definecliquesigma}). Note that, for the
first combination, $\Sigma= I_p$. Also, since the OHC is permutation
invariant, we take $\Gamma= I_p$ for simplicity so that $\Sigma=
\Sigma_0$.
For each $\Sigma$, we generate $n$ samples $X_1, X_2, \ldots, X_n$ from
$N(0, \Sigma)$ and obtain $\rho_{ij}$ for all $1 \leq i < j \leq p$.

%
\begin{figure*}

\includegraphics{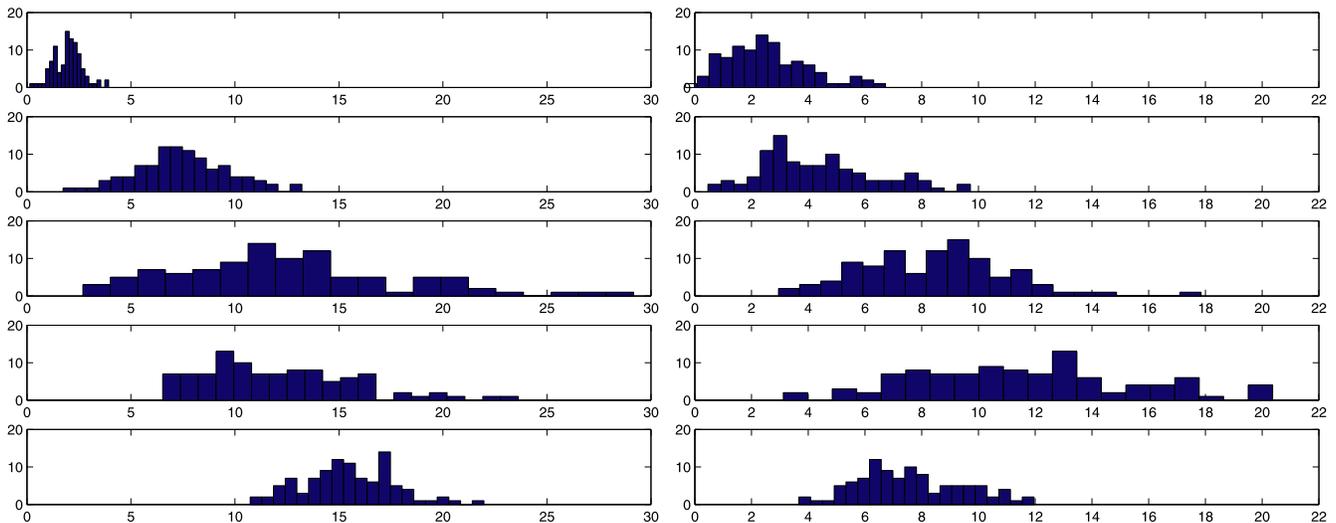}

\caption{Simulated scores of OHC applied to pairwise correlations (left
column) and maximum correlations (right column). The top panel
represents the case of no cliques. The others represent cliques with
various size and strength (introduced in the text). In comparison,
OHC applied to maximum correlations has smaller power than OHC applied
to pairwise correlations.}\label{fig:coordinatewise}
\end{figure*}

In the first approach, we sort all $N = p(p-1)/2$ different $P$-values
$\{\pi_{ij}\dvtx 1 \leq i < j \leq p\}$ in ascending order and write them
as follows:
\[
\pi_{(1)} < \pi_{(2)} < \cdots< \pi_{(N)},\quad N = p
(p-1)/2.
\]
We then apply the Orthodox Higher Criticism (OHC) and obtain the HC score
%
\begin{eqnarray}
\label{HCmatrix} \mathrm{OHC}_{N}^+& =& \max_{i\dvtx 1/N \leq\pi_{(i)} \leq1/2} \sqrt{N}
\bigl[(i/N) - \pi _{(i)} \bigr]
\nonumber
\\[-8pt]
\\[-8pt]
\nonumber
&&{}/ \sqrt{\pi_{(i)} (1 -
\pi_{(i)})}.
\end{eqnarray}
In the second approach, we sort all $p$ different $P$-values $\pi_j^*$,
$1 \leq j \leq p$, in the ascending order and denote them by
\[
\pi_{(1)}^* < \pi_{(2)}^* < \cdots< \pi_{(p)}^*.
\]
We then apply the HC by (\ref{HCmatrix}), but with $\pi_{(i)}$ replaced
by $\pi_{(i)}^*$.

The histograms of $\mathrm{OHC}_N^+$ based on $100$ repetitions are displayed in
Figure~\ref{fig:coordinatewise},
which suggests that OHC yields satisfactory detection.
For all four types of cliques, the OHC applied in the second approach
has smaller power in separation than that in the first approach.

\subsection{Detecting Low-Rank Perturbations of the Identity Matrix}
\label{sec:lowrank}
Now we test whether $\Sigma= I_p$ or instead we have a low-rank
perturbation $\Sigma= I + H$,
where the rank $r$ of $H$ is relatively small compared to $p$.\footnote
{The model $I + H$
is an instance of the so-called spiked covariance model \cite{John02};
there are of course hypothesis tests
specifically developed for this setting using random matrix theory. We
thought it
would be interesting to derive what the HC viewpoint offers in this situation.}
Consider the spectral decomposition
\[
\Sigma= Q \Lambda Q',
\]
where $Q$ is a $p\times p$ orthogonal matrix and $\Lambda$ is a
diagonal matrix, with
the first $r$ entries $1 + h_i$, $h_i > 0$, $1 \leq i \leq r$, and
other diagonal entries $1$.
We assume the eigenbasis $Q$ is unknown to us.
In a ``typical'' eigenbasis, the coordinates of $Q$ will be ``uniformly
small,'' so that even if some of the eigenvalue excesses $h_i$ are
nonzero $0$, the matrix $\Sigma$ can be
close to the corresponding coordinates of $I_p$.
Therefore, the pairwise covariances may be a very poor tool
for diagnosing departure form the null.

Instead we work with the empirical spectral decomposition
and apply HC to the sorted empirical eigenvalues.
Denote the empirical covariance matrix by
\[
S_n = (1/n) X'X,
\]
and let
\[
\lambda_{1} > \lambda_{2} > \cdots> \lambda_{n}
\]
be the (nonzero) eigenvalues of $S_n$ arranged in the descending order.
The sorted eigenvalues play a role analogous to the sorted $P$-values
in the
earlier sections, since the perturbation of $I_p$ by a low-rank matrix $H$
will inflate a small fraction of the empirical eigenvalues, similar
to the way the top few order statistics are inflated in the rare/weak model.
We define our approximate $Z$-scores by
standardizing each $\lambda_i$ using its mean and standard deviation
under the
null.
The resulting $t$-like statistics, which we call the $\mathrm{eigenHC}$, are
%
\begin{equation}
\label{DefineeHC} \mathrm{eigenHC}_{n, i} = \frac{(\lambda_i - E_0[\lambda_i])}{\mathrm
{\mathrm{SD}}_0(\lambda_i)}, \quad 1 \leq i \leq p,
\end{equation}
where $E_0[\lambda_i]$ and $\mathrm{SD}_0(\lambda_i)$ are the mean and standard
deviation of
$\lambda_i$ evaluated under the null hypothesis $\Sigma= I_p$,
respectively. Note that $E_0[\lambda_i]$ and $\mathrm{SD}_0(\lambda_i)$
can be conveniently evaluated by Monte-Carlo simulations.\footnote
{Since these are the eigenvalues of a standard Wishart matrix, much
existing analytic information is applicable.
For example, under the null distribution, the top several eigenvalues
are dependent and non-Gaussian;
Johnstone \cite{John02} showed that the distribution of the top
eigenvalue is Tracy--Widom. Here we do not use such refined mathematical
analyses, but only Monte-Carlo simulations.}

%
\begin{figure*}

\includegraphics{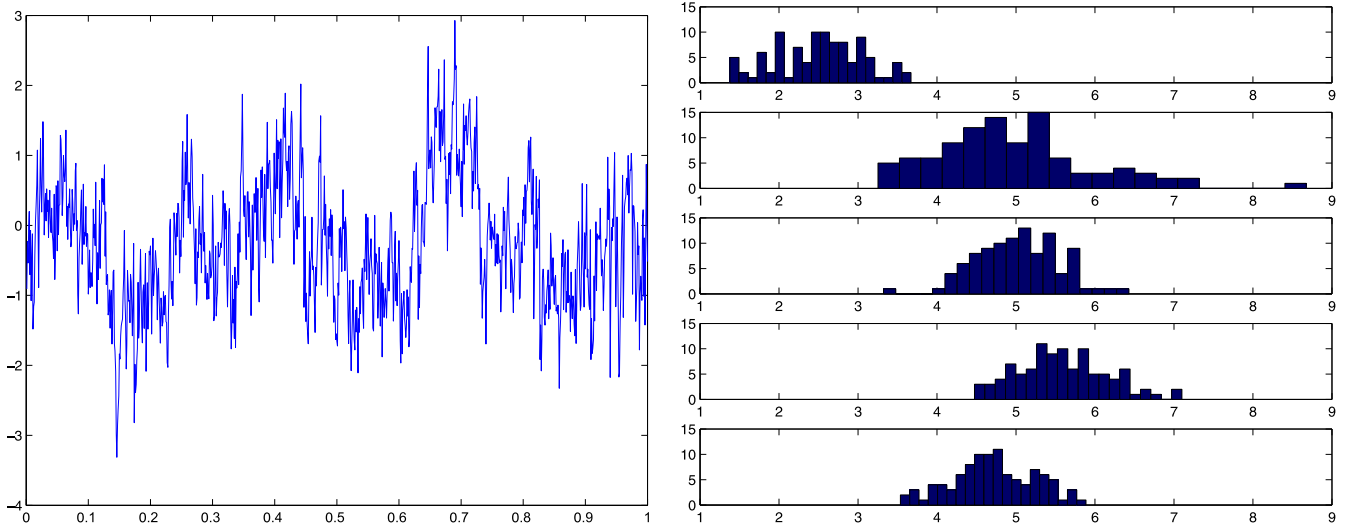}

\caption{Left column: $x$-axis is $i/n$, $y$-axis is $\mathrm{eigenHC}_{n, i}$.
Right column: simulated scores of $\mathrm{eigenHC}^*_n$.
The top panel represents the unperturbed case. The others represent
contamination with different ranks (introduced in text).}
\label{fig:eig}
\end{figure*}

In Figure~\ref{fig:eig}, we present a realization of $\{\mathrm{eigenHC}_{n,i}\dvtx\break 1\leq i\leq p\}$ in the case of $n = p = 1000$.
The figure looks vaguely similar to realizations of a normalized
uniform empirical process,
which suggests that the normalization in (\ref{DefineeHC}) makes sense.
We consider the test statistic
\[
\mathrm{eigenHC}_n^* = \max_{1 \leq i \leq\alpha_0 n} \{ \mathrm{eigenHC}_{n,i}
\},
\]
where $\alpha_0$ is a tuning parameter we set here to $1/2$.

We conducted a small-scale simulation experiment, with $(p, n) = (1000,1000)$.
For each of the $5$ different combinations of $(r,h) = (0, 0), (5, 1),\break
(15, 0.5),  (45,0.2), (135, 0.05)$, we
let $\Lambda$ be the $p \times p$ diagonal matrix with first $r$
coordinates equal to $(1 + h)$ and remaining coordinates all $1$.
We then randomly generated a
$p \times p$ orthogonal matrix $Q$ (according to the uniform measure
on orthogonal matrices) and set
\[
\Sigma= Q \Lambda Q'.
\]
Note that when\vspace*{-1pt} $(r, h) = (0,0)$, $\Sigma= I_p$.
Next, for each~$\Sigma$, we generated data $X_1, X_2, \ldots, X_n
\stackrel{\mathrm{i.i.d.}}{\sim}\break N(0,\Sigma)$ and
applied $\mathrm{eigenHC}_n^*$ to the synthetic data. Simulated results for
$100$ such synthetic data sets are reported in
Figure~\ref{fig:eig}, illustrating that HC can yield satisfactory
results even for
small $r$ or $h$.\footnote{Our point here is not that HC should
replace formal methods using random matrix
theory, but instead that HC can be used in structured settings where
theory is not yet available. A careful
comparison to formal inference using random matrix theory---not
possible here---would illustrate
the benefits of theoretical analysis of a specific situation---as
exemplified by random matrix theory, in this case---over the direct application of a general procedure like HC.}


Testing hypotheses about large covariance matrices has received much
attention in recent years.
For example, Arias-Castro {et al}. \cite{Erycluster} tested that
the underlying covariance matrix is the
identity versus the alternative where there is a small subset of
correlated components.
The correlated components may have a certain combinatorial structure
known to the statistician.
Butucea and Ingster \cite{Ingstermatrix} consider testing the null model
that the coordinates are i.i.d. $N(0,1)$ against a rare/weak model where
a small fraction of them has significantly nonzero means.
Muralidharan \cite{Muralidharan} is also related;
it adapts HC to test column dependences in gene microarray data.

%
\begin{figure*}

\includegraphics{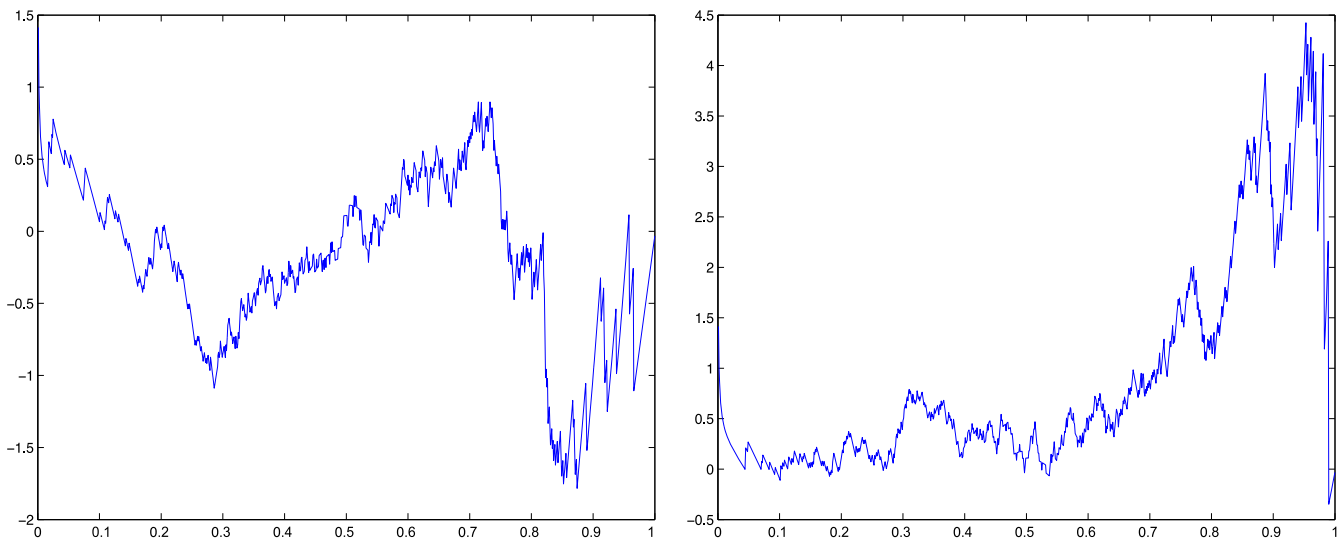}

\caption{Plot of $\mathrm{pairHC}_{n,k}$ versus $k/n$ under the null (left) and
under the alternative [right; $(\eps, \tau, \rho) = (0.05, 1, 0.25)$].
$x$-axis is $k/n$ ($n = 1000$).}
\label{fig:HCNull}
\end{figure*}

\section{Sparse Correlated Pairs Among Many Uncorrelated Pairs}

Suppose we observe independent samples $(X_i, Y_i)$, $1 \leq i \leq n$,
from a bivariate distribution with
zero means and unit variances, which is generally unknown to us.
Under the null hypothesis,
the $(X_i)$ are independent of the corresponding $(Y_i)$ (and each other);
but under the alternative, for \textit{most} pairs $(X_i, Y_i)$,
independence holds, while
for a small fraction $(X_i, Y_i)$,
the two coordinates may be correlated and each may have an elevated mean.
In short, some small collection of the pairs is correlated, unlike the
bulk of the data.

Since the underlying distribution of the pairs $(X_i, Y_i)$ is unknown
to us,
we base our test statistics on ranks $(r_i, s_i)$ of the data $(X_i, Y_i)$.
Our strategy is to compare the number of rank-pairs in the upper right
corner to the
number that would be expected under independence.

For $1 \leq k \leq n$, let
\begin{eqnarray*}
S_k& =& \#\bigl\{1 \leq i \leq n\dvtx \min\{r_i,
s_i\} \geq k\bigr\} \\
&= &\sum_{i = 1}^n
1\bigl\{ \min\{ r_i, s_i \} \geq k\bigr\}.
\end{eqnarray*}
Under the null, we have $P(\min\{r_i, s_i\} \geq k) =P(r_i \geq k)
P(s_i \geq k) = (1 - k/n)^2$, so
\[
E[S_k] = P(r_i \geq k) P(s_i \geq k) = n
(1 - k/n)^2
\]
and
\[
\operatorname{Var}(S_k) = n \bigl[(1 - k/n)^2 \bigl(1 - (1 -
k/n)^2\bigr)\bigr].
\]
Therefore, the HC idea applies as follows. Define
\[
\mathrm{pairHC}_{n,k} = \sqrt{n} \frac{S_k/n - (1 - k/n)^2}{\sqrt{(1 - k/n)^2 (1
- (1 - k/n)^2) }}
\]
and
\[
\mathrm{pairHC}_n^* = \max_{(1 - \alpha_0) n \leq k \leq n} \mathrm{pairHC}_{n,k}.
\]
Here, $\alpha_0$ is a tuning parameter and is set to $1/2$ below.


To illustrate this procedure,
suppose that $(X_i, Y_i)$, $1\leq i\leq n$, are i.i.d. samples from a
mixture of two bivariate normals
\begin{eqnarray*}
&&(1 - \eps) N(0, I_2) + \eps N(\tau{\mathbf1}_2, \Sigma),\\
&&{\mathbf 1}_2 =\pmatrix{ 1
\cr
1 }, \quad
\Sigma=\pmatrix{ 1 & \rho
\cr
\rho& 1 },
\end{eqnarray*}
where $(\eps, \tau, \rho)$ are parameters.
In Figure~\ref{fig:HCNull}, we show a plot of $\mathrm{pairHC}_{n,k}$ for $n =
1000$ and $k = 1, 2, \ldots, n$ under
the null and under the alternative where $(\eps, \tau, \rho) = (0.05,
1, 0.25)$.

%
\begin{figure*}[b]

\includegraphics{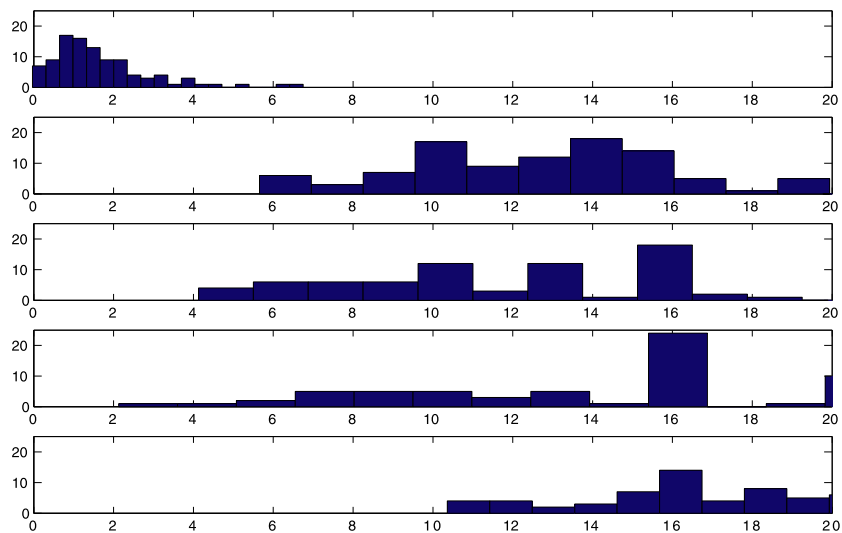}

\caption{Simulated scores of $\mathrm{pairHC}_n^*$. The top panel represents the
null case. The others represent the alternative cases with different
$(\eps, \tau, \rho)$, introduced in the text ($n = 1000$).}
\label{fig:2dhist}
\end{figure*}

We conducted a small simulation experiment as follows:
\begin{itemize}
\item Fix $n = 1000$ and define $5$ different settings where $(\eps,
\tau, \rho) = (0,0,0)$,
$(0.02,0,2.5)$, $(0.02,0.50,2)$,\break $(0.01,0.50,2.5)$ and $(0.01,0.25,3)$.
Note that the first setting corresponds to the null case.
\item Within each setting, conduct 100 Monte-Carlo repetitions, each
time generating a synthetic data set
with the given parameters and applying
$\mathrm{pairHC}_n^*$.
\end{itemize}
The results are reported in Figure~\ref{fig:2dhist},
which suggests that HC yields good separation even when the signals are
relatively rare and weak.

\section{Asymptotic Rare/Weak Model}
\label{sec:phase}

In this section we review the rare/weak signal model and discuss the
advantages of HC in this setting.

Return to the
problem (\ref{mixture1})--(\ref{mixture2}) of detecting a sparse
Gaussian mixture.
We introduce an asymptotic framework which we call the \textit{Asymptotic
Rare/Weak} (ARW) model.
We consider a sequence of problems, indexed by the number $N$
of $P$-values (or $Z$-scores, or other base statistics); in the $N$th problem,
we again consider mixtures $(1-\eps) N(0,1) + \eps N(\tau,1)$,
but now we tie the behavior of
$(\eps, \tau)$ to $N$, in order to honor the spirit of the Rare/Weak
situation.
In detail, let $\vartheta\in(0,1)$ and set
\[
\eps= \eps_N = N^{-\vartheta},
\]
so that, as $N \goto\infty$, the nonnull effects in $H_1^{(N)}$
become increasingly rare.
To counter this effect, we let $\tau_N$ tend to $\infty$ slowly, so that
the testing problem is (barely) solvable. In detail, fix $r > 0$
and set
\[
\tau_N = \sqrt{2 r \log(N)}.
\]
With these assumptions the Rare/Weak setting corresponds
to $\vartheta> 1/2$, rare enough that a shift in the overall mean is
not detectable,
and $r < 1$, weak enough that a shift in the maximum observation is not
detectable.


The key phenomenon in this model
is a \textit{threshold for detectability}. As a measure of detectability,
consider the best possible
sum of Types~I and~II errors of the optimal test.
Then, there will be a precise threshold separating values of
$(\vartheta
,r)$, where
the presence of the mixture is detectable from those where it is not
detectable.

Let
\[
\rho(\vartheta) = \cases{ %
 \vartheta- 1/2, & $1/2 <
\vartheta\leq3/4,$
\vspace*{2pt}\cr
\bigl(1 - \sqrt{(1 - \vartheta})\bigr)^2, & $3/4 < \vartheta< 1.$ }
\]
When $r > \rho(\vartheta)$,
the hypotheses separate asymptotically: the best sum of
Types~I and~II errors tends to $0$ as $N$ tends to $\infty$.
On the other hand, when
$r < \rho(\vartheta)$, the sum of Types~I and~II errors of any test
cannot get substantially smaller than $1$. The result was first proved by
Ingster \cite{Ingster97,Ingster99},
and then independently by Jin \cite{JinThesis,Jin04}.

%
\begin{figure}[t]

\includegraphics{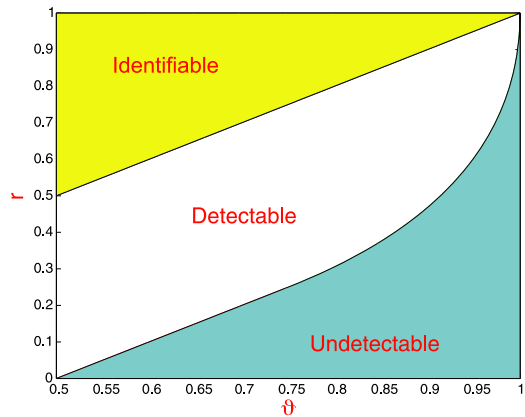}

\caption{Phase diagram for the detection problem. The detection
boundary separates the $\vartheta$-$r$ plane into the detectable region
and the undetectable region. In the identifiable region, it is not only
able to reliably tell the existence of nonzero coordinates, but also
possible to identify individual nonzero coordinates.}
\label{fig:Detect}
\end{figure}

In other words, in the two-dimensional $\vartheta$-$r$ \textit{phase space},
the curve $r = \rho(\vartheta)$ separates the bounded region
$\{(\vartheta, r)\dvtx 1/2 < \vartheta< 1, 0 < r< 1\}$ into two separate
subregions, the \textit{detectable region} and the \textit
{undetectable region}.
For $(\vartheta, r)$ in the interior of the detectable region,
two hypotheses separate asymptotically and it is possible to separate
them. For $(\vartheta, r)$ in the undetectable region,
two hypotheses merge asymptotically, and it is impossible to separate them.
Hence, the phase diagram splits into two ``phases''; see Figure~\ref{fig:Detect} for illustration.

Fix $(\vartheta, r)$ in the detectable region. Suppose we reject
$H_0^{(N)}$ if and only
\[
\mathrm{HC}_N^* \geq h(N, \alpha_N),
\]
where $\alpha_N$ tends to $0$ slowly enough so that $h(N,\break \alpha_N) =
O(\sqrt{2 \log\log(N)})$.
Then when $H_1^{(N)}$ can be detected by the optimal test, HC also
detects it, as
$N \goto\infty$.

\begin{figure*}[b]

\includegraphics{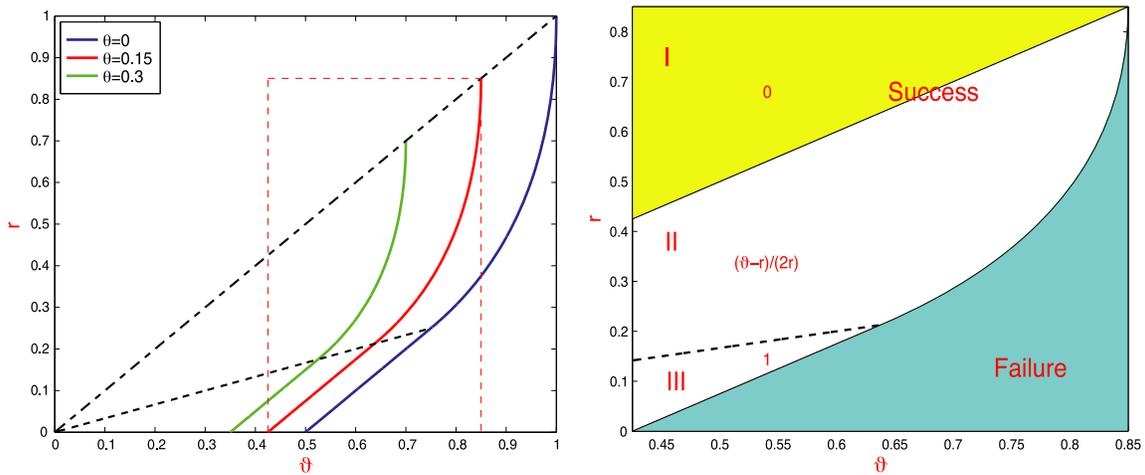}

\caption{Left: curves $r = \rho_{\theta}(\vartheta)$ for $\theta= 0,
0.15, 0.3$. Two dashed black lines are $r = \vartheta$ and $r =
\vartheta
/3$, respectively.
For $\theta= 0.15$, the most interesting region is
represented by the rectangular box.
Right: enlargement of the rectangular box. The curve $r = \rho_{\theta
}(\vartheta)$ ($\theta= 0.15$) splits the box into two subregions:
Failure (cyan) and Success (white and yellow). The two lines $r =
\vartheta$ and $r = \vartheta/3$ further split Region of Success into
three subregions, I, II and III,
where the leading terms of $q^{\mathrm{ideal}}(\vartheta, r, p)$ in (\protect\ref
{qideal}) are shown.
In the yellow region, it is not only possible to have successful
classifications, but is also possible to separate useful features from
useless ones.}
\label{fig:Classification}
\end{figure*}

Since HC can be applied without knowing the
underlying parameter $(\vartheta, r)$, we say HC is optimally adaptive.
HC thus has an advantage over the
Neyman--Pearson likelihood ratio test (LRT),
which requires precise information about the underlying ARW parameters.
The HC approach can be applied much more generally;
it is only for theoretical analysis that we focus on the narrow ARW model.

A similar phase diagram holds in
a classification problem considered in Section~\ref{sec:class},
provided that we calibrate the parameters appropriately.
Consider a sequence of classification problems indexed by $(n,p)$,
where $n$ is the number of observations and $p$ the number of
features available to the classifier.
Suppose that two classes are equally likely so that
$P(Y_i = 1) = P(Y_i = -1) = 1/2$ for all $1 \leq i \leq n$.
For $Z$ in (\ref{Zvector}), recall that $Z \sim N(\sqrt{n} \mu, I_p)$.
We calibrate with
\[
\sqrt{n} \mu(j) \stackrel{\mathrm{i.i.d.}} {\sim} (1 - \eps) \nu_0
+ \eps\nu _{\tau},
\]
where $\nu_a$ denotes the point mass at $a$.
Similarly, we use an ARW model, where we fix $(\vartheta, r, \theta)
\in(0,1)^3$ and let
\begin{eqnarray*}
\eps&= &\eps_p = p^{-\vartheta},\quad \tau= \tau_p = \sqrt{2
r \log (p)},\\
  n &=& n_p = p^{\theta}.
\end{eqnarray*}
Note that when $p \goto\infty$, $n_p$ grows with $p$, but is still
much smaller than $p$. We call such growth \textit{regular growth}. The
results below hold for other types of growth of $n$; see Jin \cite
{JinPNAS}, for example.



It turns out that there is a similar phase diagram associated with the
classification problem.
Toward this end, define
\[
\rho_{\theta}(\vartheta) = (1 - \theta) \rho\biggl(\frac{\vartheta}{1
- \theta
}
\biggr),\quad 0 < \vartheta< (1 - \theta).
\]
Fix $\theta\in(0,1)$. In the two-dimensional phase space, the most
interesting region for $(\vartheta, r)$ is the rectangular region
$\{(\vartheta, r)\dvtx 0 < \vartheta, r < (1 - \theta)\}$. The region
partitions into two subregions:
\begin{itemize}
\item(\textit{Region of success}). If $r > \rho_{\theta}(\vartheta)$,
then the HC threshold $\hat{t}_p^{\mathrm{HC}} / t_p^{\mathrm{ideal}} \goto1$ in
probability; $t_p^{\mathrm{ideal}}$ is the ideal threshold that one would choose
if the underlying parameters $(\vartheta, r)$ are known. Note that HCT
is driven by data, without the knowledge of $(\vartheta, r)$. Also, the
classification error of an HCT-trained classification rule tends to $0$
as $p \goto\infty$.
\item(\textit{Region of failure}). When $r < \rho_{\theta
}(\vartheta)$,
the classification error of any trained classification rule tends to
$1/2$, as $p \goto\infty$.
\end{itemize}
See Figure~\ref{fig:Classification}.
The above includes the case where $n_p \goto\infty$ but $n_p / p^{a}
\goto0$ for any fixed $a > 0$ as the special case of $\theta= 0$.
See more discussion in \cite{DJ08,JinPNAS,DJ09}.
Ingster {et al.} \cite{IPT} derived independently the
classification boundary, in a broader setting than that in \cite{DJ08,JinPNAS,DJ09}, but they
did not discuss HC.

The conceptual advantage of HC lies in its ability to perform optimally
under the ARW framework---without needing to know the underlying ARW parameters:
HC is a data-driven nonparametric statistic that is not tied to
the idealized model we discussed here, and yet works well in this
model.


The phase diagrams above are for settings where the test statistics or
measured features $X_i$
are independent normals with unit variances (the normal means may be different).
In more complicated settings, how to derive the phase diagrams is an
interesting but nontrivial problem.
Delaigle {et al}. \cite{Robustness} studied the problem of
detecting sparse mixtures by
extending models (\ref{mixture1})--(\ref{mixture2}) to a setting where
$X_i$ are the Student's $t$-scores based on (possibly) non-Gaussian
data, where the marginal density of $X_i$ is unknown but is
approximately normal.
Fan {et al}. \cite{FanHC} extended the classification problem
considered in Section~\ref{sec:class} to a setting where the measured
features are correlated; the covariance matrix is unknown and must be
estimated. In general, the approximation errors (either in the
underlying marginal density or the estimated covariance matrix) have a
negligible effect on the phase diagrams when the true effects are
sufficiently sparse
and a nonnegligible yet subtle effect when the true effects are
moderately sparse; see \cite{Robustness,FanHC}.

\subsection{Phase Diagram in the Nonasymptotic Context}
The phase diagrams depict an asymptotic situation;
it is natural to ask how they behave for finite $N$.
This has been studied in \cite{DJ09}, Figure~4, Sun \cite{Sun},
Blomberg \cite{Blomberg} and \cite{Rohban}.
In principle, for finite $N$, we would not experience ``perfectly
sharp'' phase transition
as visualized in Figures~\ref{fig:Detect} and \ref{fig:Classification}.
However, numeric studies reveal that, for reasonably large $N$,
the transition zone between the region where inference can be rather
satisfactory
and the region where inference is nearly impossible
is comparably narrow, increasingly so as $N$ increases.
Sun \cite{Sun} used such ideas to study a GWAS
on Parkinson's disease, and argued that standard designs for GWAS
are inefficient in many cases. Xie \cite{Xie} and Wu \cite{Meihua}
used the phase diagram as a framework for sample size and power calculations.

\subsection{Phase Diagram for FDR-Controlling Methods}
\label{subsec:phaseFDR}
Continuing the discussion in Section~\ref{subsec:featureFDR}, we
investigate the optimal FDR control parameter $q$ in the rare/weak setting.
Suppose we select features by applying Benjamini--Hochberg's
FDR-controlling method
to the ARW. The ``ideal'' FDR control parameter $q^{\mathrm{ideal}}(\vartheta,
r, p)$ is the feature-FDR associated with $t_p^{\mathrm{ideal}}$
(i.e., we have a discovery if and only if the feature $Z$-score exceeds
$t_p^{\mathrm{ideal}}$ in magnitude). In Donoho and Jin~\cite{DJ09}, it is shown
that, as $p \goto\infty$,
%
\begin{eqnarray}\quad
\label{qideal}&& q^{\mathrm{ideal}}(\vartheta, r, p)
\nonumber
\\[-8pt]
\\[-8pt]
\nonumber
&&\quad= \cases{ %
 o(1), & $r > \vartheta,$
\vspace*{2pt}\cr
\displaystyle\frac{\vartheta- r}{2r} + o(1), & $\vartheta/3 < r < \vartheta,$
\vspace*{2pt}\cr
1 - o(1), & $\rho_{\theta}(\vartheta) < r < \vartheta/3,$}
\end{eqnarray}
which gives an interesting 3-phase structure: see Regions I, II, III in
Figure~\ref{fig:Classification}.
Somewhat surprisingly, the optimal FDR is very close to $1$
in one of the three phases (i.e., Region III); in this phase, to obtain optimal
classification behavior, we
set the feature selection threshold very low so that we include most of the
useful features; but when we do this, we necessarily include many
useless features, which dominate in numbers among all selected features.
Similar comments \mbox{apply} when replacing Benjamini--Hochberg's FDR control
by the local FDR (Lfdr) approach of Efron \cite{Efron}; see \cite{DJ09}
for details.


%

\subsection{Phase Diagrams in Other Rare/Weak Settings}
Phase diagrams offer a new criterion for measuring performance in
multiple testing in the rare/weak effects model.

This framework is useful in many other settings.
Consider a linear regression model $Y = X \beta+ z$, $z \sim N(0,
I_n)$, where $X = X_{n,p}$ and $p \geq n$.
The signals in the coefficient vector $\beta$ are rare and weak,
and the goal is variable selection (different from that in Section~\ref{subsec:corr}).
In a series of papers \cite{GJW,JiJin,JZZ,Ke},
we use the Hamming error as the loss function for variable selection
and study phase diagrams
in settings where the matrices $X$ get increasingly ``bad'' so the
problem get increasingly harder.
These studies propose several new variable selection procedures,
including UPS \cite{JiJin}, Graphlet Screening \cite{JZZ} and CASE
\cite{Ke}.
The study is closely related to \cite{Donoho} on Compressed Sensing.

References \cite{JW1,JW2} present ARW phase diagrams for sparse spectral clustering,
and \cite{Fienberg} presents phase diagrams for computer privacy and
confidentiality.

\section*{Acknowledgments}
The authors thank the Associate Editor and referees for very helpful comments.

To the memory of John W. Tukey 1915--2000 and of Yuri I. Ingster 1946--2012,
two pioneers in mathematical statistics. The research of D. Donoho
supported in part by NSF Grant DMS-09-06812 (ARRA),
and the research of J. Jin supported in part by NSF Grant DMS-12-08315.


%
%

%


\end{document}